\documentclass[12pt, reqno]{amsart}
\usepackage{amsmath, amsthm, amssymb, mathtools}
\usepackage{enumitem}
\usepackage{cite}
\usepackage[margin=2.5cm]{geometry}
\usepackage{xcolor}
\usepackage[colorlinks=true, citecolor=blue, linkcolor=blue, urlcolor=blue]{hyperref}

\newtheorem{theorem}{Theorem}[section]
\newtheorem{lemma}[theorem]{Lemma}
\newtheorem{proposition}[theorem]{Proposition}

\theoremstyle{definition}
\newtheorem{definition}[theorem]{Definition}
\theoremstyle{remark}
\newtheorem{remark}[theorem]{Remark}
\numberwithin{equation}{section}

\newcommand{\R}{\mathbb{R}}
\newcommand{\N}{\mathbb{N}}
\newcommand{\RN}{\mathbb{R}^N}
\newcommand{\cL}{\mathcal{L}}
\newcommand{\cE}{\mathcal{E}}
\newcommand{\cB}{\mathcal{B}}
\newcommand{\cS}{\mathcal{S}}
\newcommand{\LN}{\Lambda_N}
\newcommand{\Cs}{C_{N,s}}
\newcommand{\dd}{\,\mathrm{d}}
\DeclareMathOperator{\supp}{supp}
\DeclareMathOperator*{\essinf}{ess\,inf}

\begin{document}

\title[Singular local--nonlocal problems with Hardy potential]{Mixed local--nonlocal problems with Hardy potential and singular nonlinearity}

\author[S. Ghosh]{Sekhar Ghosh}

\address{Sekhar Ghosh: Department of Mathematics, \newline National Institute of Technology Calicut, Kozhikode, Kerala, 673601, India}
\email{sekharghosh1234@gmail.com / sekharghosh@nitc.ac.in}

\author[V. Kumar]{Vishvesh Kumar}
\address{Vishvesh Kumar: Department of Mathematical Sciences,\newline Indian Institute of Technology (BHU) Varanasi, Varanasi, Uttar Pradesh, 221005, India}
\email{vishvesh.mat@iitbhu.ac.in/ vishveshmishra@gmail.com}

\author[A. P. Verma]{Amar Pal Verma}
\address{Amar Pal Verma: Department of Mathematics,\newline National Institute of Technology Rourkela, Rourkela, Odisha, 769008, India}
\email{amarpalverma100@gmail.com}

\date{\today}

\begin{abstract}
We study the existence and the regularity of positive solutions of the singular problem
\[
-\Delta u+(-\Delta)^s u-\mu\frac{u}{|x|^2}=\frac{g(x)}{u^{\gamma}}\ \text{ in }\Omega,\qquad u>0\ \text{ in }\Omega,\qquad u=0\ \text{ in }\RN\setminus\Omega,
\]
where $\Omega\subset\RN$, $N\ge 3$, is a bounded domain containing the origin, $s\in(0,1)$, $\gamma\ge 0$, $0<\mu<\LN=(N-2)^2/4$, and $g\ge 0$ belongs to $L^m(\Omega)$ for some $m\ge1$. Depending on the values of $m$ and $\gamma,$ we prove the existence of a positive weak solution which belongs to $W^{1,q}_0(\Omega)$ for an explicit $q<2$, to the energy space $H^1_0(\Omega)$, or to $H^1_{\rm loc}(\Omega)$ with $u^{(\gamma+1)/2}\in H^1_0(\Omega)$. In the infinite-energy regimes, the admissible size of the Hardy parameter is given by an explicit threshold depending on $N$, $m$ and $\gamma.$ In particular, for $\gamma=0,$ it coincides with the threshold known in the summability theory of the purely local problem. We also prove an optimal-type summability result $u\in L^{\frac{Nm(1+\gamma)}{N-2m}}(\Omega)$ for $m<N/2$ and $\mu$ below an explicit threshold $\mu_{m,\gamma}$, uniqueness of finite-energy solutions, and nonexistence of locally finite-energy solutions for $\mu>\LN$, the latter by means of a Picone-type inequality for the mixed operator. Moreover, we show that every solution satisfies $u\ge A|x|^{-\alpha_\mu}$ near the origin, with $\alpha_\mu=\frac{N-2}{2}-\sqrt{\LN-\mu}$, which proves that the summability threshold is sharp in a range of parameters, and we prove existence in the critical case $\mu=\LN$ by means of improved Hardy inequalities. The existence proofs combine an approximation scheme, a priori estimates obtained by testing with suitable powers of the approximate solutions together with the Hardy inequality, and a weak Harnack inequality for the mixed operator, which provides positivity bounds independent of the approximation parameter.
\end{abstract}

\maketitle
\allowdisplaybreaks

\section{Introduction}\label{sec:intro}
Let $\Omega\subset\RN$, $N\ge3$, be a bounded domain (that is,  a bounded, open, and connected set) with $0 \in \Omega$, let $s\in(0,1)$, $\gamma\ge 0$, $\mu>0$, and let $\LN:=\frac{(N-2)^2}{4}$
be the Hardy constant. In this paper, we study the singular problem
\begin{equation}\tag{$\mathcal{P}_\mu$}\label{problem1}
\left\{
\begin{aligned}
\cL u-\mu\frac{u}{|x|^{2}}&=\frac{g(x)}{u^{\gamma}} &&\text{in } \Omega,\\
u&>0 &&\text{in } \Omega,\\
u&=0 &&\text{in } \RN\setminus\Omega,
\end{aligned}
\right.
\end{equation}
where $g\in L^m(\Omega)$, $m\ge1$, is nonnegative and not identically zero, and $\cL$ is the mixed local--nonlocal operator
\[
\cL:=-\Delta+(-\Delta)^s,\qquad
(-\Delta)^s u(x)=\Cs\,\lim_{\varepsilon\to0^+}\int_{\{|x-y|\ge\varepsilon\}}\frac{u(x)-u(y)}{|x-y|^{N+2s}}\dd y ,
\]
$\Cs>0$ being the usual normalization constant of the fractional Laplacian. Our existence results concern $0<\mu\le\LN$; the main case is $0<\mu<\LN$, the critical case $\mu=\LN$ is treated in Theorem \ref{thm:critical}, and nonexistence for $\mu>\LN$ is proved in Theorem \ref{thm:nonexistence}.

Operators of the form $\cL$ are a particular instance of a class of integro-differential operators which has been studied since the late 1960s in connection with the validity of the maximum principle \cite{bony1968semigroupes,cancelier1986problemes}; from the probabilistic point of view, $\cL$ is the generator of the superposition of a Brownian motion and an independent $2s$-stable L\'evy process, see e.g. \cite{chen2012boundary}. Such superpositions appear in models in which a population, or a diffusing particle, may move both by local random walks and by long jumps, see e.g. \cite{dipierro2022nonlocal} and the references therein. In recent years a systematic analytical theory for $\cL$ has been developed, including existence results, maximum principles, interior and boundary regularity, and Harnack-type inequalities; see e.g. \cite{antonini2025global, biagi2022mixed, biagi2023faber, defilippis2024gradient, garain2022regularity, garain2023higher, su2022regularity, su2025some}.

Elliptic problems with singular nonlinearities, that is, with a right-hand side which blows up where the solution vanishes, have a long history, starting with the seminal papers \cite{crandall1977dirichlet, lazer1991singular}. We recall, in particular, that, for smooth data bounded away from zero, the solution of $-\Delta u=u^{-\gamma}$ belongs to $H^1_0(\Omega)$ if and only if $\gamma<3$ \cite{lazer1991singular}. For data in Lebesgue spaces $L^m(\Omega)$, a complete picture of the existence and of the optimal regularity of solutions of $-\Delta u=g\,u^{-\gamma}$ in terms of $m$ and $\gamma$ was obtained by Boccardo and Orsina \cite{boccardo2010semilinear}. The analogous questions for the fractional Laplacian were addressed in \cite{barrios2015semilinear, canino2017nonlocal}, see also \cite{arora2022large}, and, for the mixed operator $\cL$, in \cite{arora2023combined, biagi2024multiplicity, garain2023class, BG2024}.

The inverse-square potential $\mu|x|^{-2}$ has the same homogeneity as the Laplacian and therefore it cannot be treated as a lower-order perturbation of $-\Delta$. Its effect is governed by the famous Hardy inequality
\begin{equation}\label{eq:hardy-intro}
\LN\int_{\RN}\frac{u^2}{|x|^2}\dd x\le\int_{\RN}|\nabla u|^2\dd x,\qquad u\in C^\infty_c(\RN),
\end{equation}
whose constant $\LN$ is optimal and not attained \cite{hardy1988inequalities, peral2021elliptic}. For $\mu>\LN$ there are no positive (even very weak) solutions of $-\Delta u-\mu|x|^{-2}u=f\ge0$ near the origin, as first observed by Baras and Goldstein \cite{baras1984heat} in their study of the associated heat equation (see also \cite{peral2021elliptic} for the stationary problem), whereas for $\mu\le\LN$ solutions exist but are in general singular at the origin \cite{brezis2005semilinear, dupaigne2002nonlinear}, and the Hardy term deeply affects their summability \cite{boccardo2006existence, brezis1997blow, garcia1998hardy, vazquez2000hardy}. In particular, Boccardo, Orsina and Peral \cite{boccardo2006existence} showed that, for $f\in L^m(\Omega)$ with $1<m<N/2$, the classical gain of summability $u\in L^{m^{**}}(\Omega)$ with $m^{**}:=\frac{Nm}{N-2m}, $ survives only if $\mu$ is smaller than an explicit threshold depending on $m$. Hardy-type problems for the fractional Laplacian were studied, among many others, in \cite{abdellaoui2016hardy, frank2008hardy}. For the mixed operator, Hardy inequalities and linear and semilinear problems with Hardy potential were recently investigated in \cite{biagi2025mixed}; see also \cite{aikyn2026mixed, biagi2022brezis, biagi2024brezis, malhotra2026eigenvalues} for related critical and sublinear problems driven by $\cL$.

To the best of our knowledge, problem \eqref{problem1}, in which the singular nonlinearity, the inverse-square potential and the mixed operator interact, has not been considered so far. Our goal is to extend the theory of \cite{boccardo2010semilinear} to this setting, and to quantify how the Hardy potential restricts the admissible range of the parameter $\mu$ in terms of $m$ and $\gamma$.

Throughout the paper, we use the following notation.
\begin{itemize}[leftmargin=*]
\item For $1<r<\infty$ we write $r'=r/(r-1)$, with $1'=\infty$ and $\infty'=1$; $2^*=2N/(N-2)$ is the critical Sobolev exponent.
\item $\cS:=\{(m,\gamma):\ 1\le m\le\infty,\ \gamma\ge0,\ (m,\gamma)\neq(1,0)\}$.
\item For $0\le\gamma<1$ we set
\begin{equation}\label{eq:msharp}
m^{\#}:=\Big(\frac{2^*}{1-\gamma}\Big)'=\frac{2N}{N+2+\gamma(N-2)}\in\Big(1,\,\frac{2N}{N+2}\Big] .
\end{equation}
\item For $(m,\gamma)\in\cS$ with $1\le m<N/2$ we set
\begin{equation}\label{eq:exponents}
\delta_m:=\frac{N(m-1)+\gamma m(N-2)}{N-2m},\qquad
\sigma_m:=\frac{Nm(1+\gamma)}{N-2m},
\end{equation}
and
\begin{equation}\label{eq:threshold}
\mu_{m,\gamma}:=\frac{4\delta_m}{(1+\delta_m)^2}\,\LN=\frac{\big[N(m-1)+\gamma m(N-2)\big](N-2m)}{m^2(1+\gamma)^2}.
\end{equation}
\item If, moreover, $0\le\gamma<1$ and $1\le m<m^\#$, we set
\begin{equation}\label{eq:qm}
q_m:=\frac{Nm(1+\gamma)}{N-m(1-\gamma)}.
\end{equation}
\end{itemize}
Elementary properties of these quantities are collected in Lemma \ref{lem:exponents} below. In particular,  $\delta_m\ge\gamma$, $\delta_m>0$, $0<\mu_{m,\gamma}\le\LN$, and, for $0\le\gamma<1$ and $m<m^\#$, $\delta_m<1$, $\sigma_m=q_m^*$, $\frac{N}{N-1}<q_m<2$ and $\mu_{m,\gamma}<\LN$. Note that $m^{\#}=(2^*)'$ when $\gamma=0$; since $(1,0)\notin\cS$, in this case Theorem \ref{thm:main1} below concerns the range $1<m<(2^*)'$.

We use the bilinear form associated with $\cL$,
\begin{equation}\label{eq:form}
\cE(u,\psi):=\int_{\RN}\nabla u\cdot\nabla\psi\dd x+\frac{\Cs}{2}\iint_{\R^{2N}}\frac{(u(x)-u(y))(\psi(x)-\psi(y))}{|x-y|^{N+2s}}\dd x\dd y,
\end{equation}
which is well defined, for instance, when $u,\psi\in H^1_0(\Omega)$ (extended by zero outside $\Omega$); see Section \ref{sec:prelim} for more details.

\begin{definition}\label{def:weak}
Let $g\in L^1(\Omega)$, $g\ge0$. A measurable function $u:\RN\to\R$ is a \emph{weak solution} of \eqref{problem1} if
\begin{enumerate}[label=(\alph*)]
\item $u=0$ a.e. in $\RN\setminus\Omega$ and $u\in L^1(\Omega)\cap W^{1,q}_{\rm loc}(\Omega)$ for some $q>\frac{N}{N-1}$;
\item for every $\omega\Subset\Omega$ there exists $c_\omega>0$ such that $u\ge c_\omega$ a.e. in $\omega$;
\item for every $\psi\in C^\infty_c(\Omega)$,
\begin{equation}\label{eq:weakform}
\cE(u,\psi)-\mu\int_\Omega\frac{u\,\psi}{|x|^2}\dd x=\int_\Omega\frac{g\,\psi}{u^\gamma}\dd x .
\end{equation}
\end{enumerate}
\end{definition}

Since $\Omega$ is bounded, $W^{1,q}_{\rm loc}(\Omega)\subset W^{1,r}_{\rm loc}(\Omega)$ for $r<q$, so in (a) we may and do assume $q\in(\frac{N}{N-1},N)$. By Lemma \ref{lem:welldef} and Remark \ref{rem:welldef} below, all the integrals in \eqref{eq:weakform} are absolutely convergent under (a) and (b). The lower bound $q>\frac{N}{N-1}$ is needed only for the Hardy term, and it cannot be relaxed to $q=\frac{N}{N-1}$: for $\frac{N-1}{N}<b\le1$, the function $|x|^{2-N}\big(\log\frac1{|x|}\big)^{-b}$ belongs to $W^{1,N/(N-1)}$ near the origin, but $u|x|^{-2}$ is not integrable there. Such fundamental-solution-type behaviour occurs for data in $L^1(\Omega)$ when $\gamma=0$, which is why the pair $(m,\gamma)=(1,0)$ is excluded from $\cS$; for $\gamma>0$ the singular term prevents it, see the exponent $q_1=\frac{N(1+\gamma)}{N-1+\gamma}>\frac{N}{N-1}$ in \cite{boccardo2010semilinear}. As in \cite{boccardo2010semilinear}, the boundary condition is encoded in the global regularity of $u$ given by the following theorems.

\begin{theorem}[Infinite energy, $0\le\gamma<1$]\label{thm:main1}
Let $0\le\gamma<1$, $(m,\gamma)\in\cS$ with $1\le m<m^\#$, and let $g\in L^m(\Omega)$ be nonnegative and not identically zero. If
\[
0<\mu<\mu_{m,\gamma},
\]
then \eqref{problem1} has a weak solution $u$ such that
\[
u\in W^{1,q_m}_0(\Omega),\qquad u^{\frac{1+\delta_m}{2}}\in H^1_0(\Omega),\qquad u\in L^{\sigma_m}(\Omega).
\]
\end{theorem}

\begin{theorem}[Finite energy]\label{thm:main2}
Let $(m,\gamma)\in\cS$ with either $\gamma=1$, or $0\le\gamma<1$ and $m\ge m^\#$. Let $g\in L^m(\Omega)$ be nonnegative and not identically zero, and let $0<\mu<\LN$. Then \eqref{problem1} has a weak solution $u\in H^1_0(\Omega)$.
\end{theorem}

\begin{theorem}[Strongly singular case, $\gamma>1$]\label{thm:main3}
Let $\gamma>1$ and let $g\in L^1(\Omega)$ be nonnegative and not identically zero. If
\[
0<\mu<\frac{4\gamma}{(1+\gamma)^2}\,\LN,
\]
then \eqref{problem1} has a weak solution $u$ such that
\[
u\in H^1_{\rm loc}(\Omega),\qquad u^{\frac{\gamma+1}{2}}\in H^1_0(\Omega).
\]
\end{theorem}

Theorems \ref{thm:main1}--\ref{thm:main3} cover every $(m,\gamma)\in\cS$, and the regularity they provide is the one obtained in \cite{boccardo2010semilinear} for $\mu=0$ and $\cL=-\Delta$.

\begin{remark}\label{rem:thresholds}
\begin{enumerate}[label=(\roman*), leftmargin=*]
\item Since $t\mapsto 4t/(1+t)^2$ is increasing on $(0,1)$ and $m\mapsto\delta_m$ is increasing, for $0\le\gamma<1$ the threshold $\mu_{m,\gamma}$ increases with $m\in[1,m^\#)$; moreover
\[
\mu_{1,\gamma}=\frac{4\gamma}{(1+\gamma)^2}\,\LN\quad(\gamma>0),\qquad \lim_{m\to m^\#}\mu_{m,\gamma}=\LN .
\]
Hence, the restriction on $\mu$ disappears continuously as one approaches the finite-energy regime of Theorem \ref{thm:main2}. The threshold in Theorem \ref{thm:main3} is invariant under $\gamma\mapsto1/\gamma$ and therefore matches $\mu_{1,1/\gamma}$.
\item For $\gamma=0$, \eqref{eq:threshold} reduces to $\mu_{m,0}=N(m-1)(N-2m)/m^2$, which is the threshold appearing in the summability theory for $-\Delta u-\mu|x|^{-2}u=f\in L^m(\Omega)$, see \cite{boccardo2006existence}.
\item The thresholds come from the fact that a priori estimates are obtained by testing with powers $u_n^\theta$, $\theta\neq1$, of approximate solutions. In this case, the Hardy term must be absorbed by $\frac{4\theta}{(1+\theta)^2}\|\nabla u_n^{(\theta+1)/2}\|_2^2$ rather than by $\|\nabla u_n\|_2^2$. Note that the nonlocal part of $\cL$ does not improve the constant: the largest constant $C$ such that $C\int_\Omega u^2|x|^{-2}\dd x\le\cE(u,u)$ for all $u\in H^1_0(\Omega)$ is still $\LN$ (see Remark \ref{rem:hardy-optimal}). For $\delta_m\ge1$ the threshold $\mu_{m,\gamma}$ is sharp for the summability of solutions (Theorem \ref{thm:origin}(c)); for $\delta_m<1$ its sharpness for \eqref{problem1} is an interesting open question.
\end{enumerate}
\end{remark}

Let us now describe the method that we follow. We use the approximation scheme of \cite{boccardo2010semilinear}. In this technique,  we truncate $g$, regularize the singular term and the potential, and consider the problems
\begin{equation}\tag{$P_n$}\label{problem11}
\left\{
\begin{aligned}
\cL u_n-\mu\frac{u_n}{|x|^{2}+1/n}&=\frac{g_n}{(u_n+1/n)^\gamma} &&\text{in } \Omega,\\
u_n&=0 &&\text{in } \RN\setminus\Omega,
\end{aligned}
\right.
\qquad g_n:=\min\{g,n\}.
\end{equation}
The main steps of our analysis to prove the existence results are the following.
\begin{enumerate}[label=(\arabic*), leftmargin=*]
\item Problem \eqref{problem11} has a unique nonnegative bounded weak solution (Lemma \ref{lem:Pn}); this follows from the Lax--Milgram theorem, the Hardy inequality and the Schauder fixed point theorem, while boundedness follows from Stampacchia-type estimates for $\cL$ (Lemma \ref{lem:stampacchia}).
\item The sequence $(u_n)$ is nondecreasing and, by a weak Harnack inequality for $\cL$, it is bounded from below by a positive constant on every compact subset of $\Omega$, uniformly in $n$ (Lemmas \ref{lem:strongmin} and \ref{lem:monotone}).
\item Testing \eqref{problem11} with suitable powers of $u_n$, and controlling the nonlocal term by an elementary algebraic inequality (Lemma \ref{lem:alg}) and the Hardy term by \eqref{eq:hardy-intro}, we obtain uniform estimates in the spaces appearing in Theorems \ref{thm:main1}--\ref{thm:main3} (Section \ref{sec:estimates}).
\item We pass to the limit in the weak formulation. For the nonlocal term this requires some care, since in Theorems \ref{thm:main1} and \ref{thm:main3} the solutions do not have finite energy. This is handled in Lemma \ref{lem:welldef}.
\end{enumerate}
The nonexistence result relies on a Picone-type inequality for $\cL$ (Lemma \ref{lem:picone}), which, combined with the optimality of $\LN$ for the mixed energy (Remark \ref{rem:hardy-optimal}), shows that a positive supersolution forces the Hardy inequality with constant $\mu$.

 The approximation scheme \eqref{problem11}, the monotonicity of $(u_n)$, the choice of the test functions $u_n^{\delta}$ and of the exponents $\delta_m$, $\sigma_m$, $q_m$, and the H\"older argument leading to the $W^{1,q}_0$ estimate are those introduced by Boccardo and Orsina \cite{boccardo2010semilinear} for $-\Delta u=g\,u^{-\gamma}$; their adaptation to the mixed operator without Hardy potential is due to Garain \cite{garain2023class} and Arora and R\u{a}dulescu \cite{arora2023combined}. The interplay between powers of the solution and the Hardy inequality, which produces the thresholds $\mu_{m,\gamma}$, is modelled on the work of Boccardo, Orsina and Peral \cite{boccardo2006existence} for the case $\gamma=0$ and $\cL=-\Delta$.

We complement the existence results above with the theorems concerning the summability of the solutions, uniqueness in the energy space, and the sharpness of the condition $\mu<\LN$.

The first one is the analogue of the summability results of \cite{boccardo2010semilinear} for $m<N/2$, and, for $\gamma=0$, of the results of \cite{boccardo2006existence} for the purely local problem with Hardy potential. Note that for $\delta_m\ge1$ (that is, for $\gamma\ge1$, or for $0\le\gamma<1$ and $m\ge m^\#$) the threshold $\mu_{m,\gamma}$ is \emph{decreasing} in $m$, and $\mu_{m,\gamma}\to0$ as $m\to N/2$.

\begin{theorem}[Summability]\label{thm:summability} Let $(m,\gamma)\in\cS$ with $1\le m<N/2$, let $g\in L^m(\Omega)$ be nonnegative and not identically zero, and assume that $ 0<\mu<\mu_{m,\gamma}. $ Then problem~\eqref{problem1} admits a weak solution $u$ such that $ u^{\frac{1+\delta_m}{2}}\in H^1_0(\Omega) $ and $$ u\in L^{\sigma_m}(\Omega), \qquad \sigma_m:=\frac{Nm(1+\gamma)}{N-2m}. $$ Moreover, according to the values of $m$ and $\gamma$, the solution $u$ coincides with the one constructed in Theorem~\ref{thm:main1}, Theorem~\ref{thm:main2}, or Theorem~\ref{thm:main3}. \end{theorem}

For $0\le\gamma<1$ and $m^\#<m<N/2$, Theorem~\ref{thm:summability} improves the integrability that follows from the conclusion $u\in H^1_0(\Omega)$ of Theorem~\ref{thm:main2}. Indeed, the Sobolev embedding yields $ u\in L^{2^*}(\Omega), $ whereas Theorem~\ref{thm:summability} gives $ u\in L^{\sigma_m}(\Omega), \,\, \sigma_m>2^*, $ because $\delta_m>1$ in this range. Similarly, when $\gamma>1$ and $m>1$, it improves the summability conclusion of Theorem~\ref{thm:main3}, since $ \sigma_m>\sigma_1 = \frac{(1+\gamma)2^*}{2}. $ In the nonsingular case $\gamma=0$, the theorem gives $ u\in L^{m^{**}}(\Omega)$ with $m^{**}:=\frac{Nm}{N-2m}, $ under the condition $ \mu< \frac{N(m-1)(N-2m)}{m^2}, $ which is the same condition as in the purely local case; see \cite{boccardo2006existence}. Finally, the fact that $ \mu_{m,\gamma}\longrightarrow0 $ as  $m\to\frac N2 $ is consistent with the well-known phenomenon that, in the presence of a Hardy potential, solutions may be unbounded near the origin even for smooth data; see, for example, \cite{boccardo2006existence,peral2021elliptic} for the local case.

The second result shows that the finite-energy solution of Theorem \ref{thm:main2} is unique and satisfies the equation for all test functions in the energy space.

\begin{theorem}[Energy formulation and uniqueness]\label{thm:uniqueness}
Let $\gamma\ge0$, let $g\in L^1(\Omega)$ be nonnegative and not identically zero, and let $0<\mu<\LN$. If $u\in H^1_0(\Omega)$ is a weak solution of \eqref{problem1}, then $g\,\psi\,u^{-\gamma}\in L^1(\Omega)$ and \eqref{eq:weakform} holds for every $\psi\in H^1_0(\Omega)$. Moreover, \eqref{problem1} has at most one weak solution in $H^1_0(\Omega)$. In particular, the solution given by Theorem \ref{thm:main2} is the unique weak solution of \eqref{problem1} belonging to $H^1_0(\Omega)$.
\end{theorem}

Finally, we show that the restriction $\mu<\LN$ cannot be dropped, at least in the class of solutions which are locally of finite energy (this class contains the solutions of Theorems \ref{thm:main2} and \ref{thm:main3}). This is the mixed local--nonlocal counterpart of the classical nonexistence result of Baras and Goldstein \cite{baras1984heat}.

\begin{theorem}[Nonexistence for $\mu>\LN$]\label{thm:nonexistence}
Let $\gamma\ge0$, let $g\in L^1(\Omega)$ be nonnegative, and let $\mu>\LN$. Then there is no function $u\in L^1(\Omega)\cap H^1_{\rm loc}(\Omega)$, with $u=0$ a.e. in $\RN\setminus\Omega$, satisfying condition (b) of Definition \ref{def:weak} and
\begin{equation}\label{eq:supersol}
\cE(u,\psi)-\mu\int_\Omega\frac{u\,\psi}{|x|^2}\dd x\ge0\qquad\text{for every }\psi\in C^\infty_c(\Omega),\ \psi\ge0 .
\end{equation}
In particular, \eqref{problem1} has no weak solution $u\in H^1_{\rm loc}(\Omega)$.
\end{theorem}

For $0<\mu\le\LN$ let
\begin{equation}\label{eq:alphamu}
\alpha_\mu:=\frac{N-2}{2}-\sqrt{\LN-\mu}\in\Big(0,\frac{N-2}{2}\Big],
\end{equation}
which is the smaller root of $t(N-2-t)=\mu$; the function $|x|^{-\alpha_\mu}$ satisfies $-\Delta|x|^{-\alpha_\mu}-\mu|x|^{-\alpha_\mu-2}=0$ in $\RN\setminus\{0\}$. In the purely local case it is well known that positive supersolutions of $-\Delta u-\mu|x|^{-2}u\ge0$ are bounded from below by a multiple of $|x|^{-\alpha_\mu}$ near the origin, we refer to \cite{peral2021elliptic}. We show that the same holds for the mixed operator, the nonlocal part being of lower order near the origin, and that this lower bound makes the thresholds of Theorem \ref{thm:summability} optimal whenever $\delta_m\ge1$.

\begin{theorem}[Behaviour near the origin]\label{thm:origin}
Let $\gamma\ge0$, let $g\in L^1(\Omega)$ be nonnegative and not identically zero, and let $0<\mu<\LN$.
\begin{enumerate}[label=(\alph*)]
\item If $u\in L^1(\Omega)\cap H^1_{\rm loc}(\Omega)$, with $u=0$ a.e. in $\RN\setminus\Omega$, satisfies (b) of Definition \ref{def:weak} and \eqref{eq:supersol}, then there exist $R>0$ and $A>0$ such that
\[
u(x)\ge A\,|x|^{-\alpha_\mu}\qquad\text{for a.e. }x\in B_R .
\]
\item The same conclusion holds for the solutions given by Theorems \ref{thm:main1}, \ref{thm:main2} and \ref{thm:main3}. In particular, none of them belongs to $L^p(B_r)$ if $r>0$, $B_r\subset\Omega$, and $p\ge N/\alpha_\mu$; they are unbounded near the origin, whatever the summability of $g$.
\item \emph{(Optimality of Theorem \ref{thm:summability}.)} Let $(m,\gamma)\in\cS$ with $1\le m<N/2$ and $\delta_m\ge1$, let $g\in L^m(\Omega)$ be nonnegative and not identically zero, and let
\[
\mu_{m,\gamma}\le\mu<\LN,\qquad\text{and, if }\gamma>1,\quad \mu<\frac{4\gamma}{(1+\gamma)^2}\LN .
\]
Then the solution of \eqref{problem1} given by Theorem \ref{thm:main2} (if $\gamma\le1$) or Theorem \ref{thm:main3} (if $\gamma>1$) does not belong to $L^{\sigma_m}(\Omega)$.
\end{enumerate}
\end{theorem}

Thus, for $\gamma\ge1$, and for $0\le\gamma<1$ with $m^\#\le m<N/2$, the threshold $\mu_{m,\gamma}$ in Theorem \ref{thm:summability} is sharp. For $0\le\gamma<1$ and $m<m^\#$ the lower bound is compatible with Theorem \ref{thm:main1} (indeed, since $\delta_m<1$, we have $\alpha_\mu<\frac{N-2}{2}<\frac{N-2}{1+\delta_m}=\frac{N}{\sigma_m}$ for every $\mu<\LN$), and the sharpness of $\mu_{m,\gamma}$ remains open.

Our last result concerns the critical value $\mu=\LN$. In this case the quadratic form $\cE(u,u)-\LN\int_\Omega u^2|x|^{-2}\dd x$ is no longer coercive on $H^1_0(\Omega)$, but, by the improved Hardy inequalities of Brezis--V\'azquez \cite{brezis1997blow} and V\'azquez--Zuazua \cite{vazquez2000hardy}, it still controls the $W^{1,q}_0(\Omega)$ norm for every $q<2$. This allows us to treat the finite-energy regime of Theorem \ref{thm:main2}, except the endpoint $m=m^\#$.

\begin{theorem}[Critical case $\mu=\LN$]\label{thm:critical}
Let $\mu=\LN$ and let $(m,\gamma)\in\cS$ with either $\gamma=1$, or $0\le\gamma<1$ and $m>m^\#$. Let $g\in L^m(\Omega)$ be nonnegative and not identically zero. Then \eqref{problem1} has a weak solution $u$ such that
\begin{enumerate}[label=(\roman*)]
\item $u\in W^{1,q}_0(\Omega)$ for every $q\in[1,2)$, $[u]_s<\infty$, and $u\in L^p(\Omega)$ for every $p<2^*$;
\item for every $\beta<\frac{N-2}{2}$ there exist $R,A>0$ such that $u\ge A|x|^{-\beta}$ a.e. in $B_R$. In particular, $u\notin L^p(B_r)$ for every $r>0$ with $B_r\subset\Omega$ and every $p>2^*$.
\end{enumerate}
\end{theorem}

Theorems \ref{thm:main2}, \ref{thm:nonexistence} and \ref{thm:critical} show that, when $\gamma=1$ or $\gamma<1$ and $m>m^\#$, solutions exist for $0<\mu\le\LN$, whereas for $\mu>\LN$ there are no solutions which are locally of finite energy. The endpoint $m=m^\#$ for $\mu=\LN$, the cases $\gamma>1$ or $m<m^\#$ for $\mu$ above the respective thresholds, and the nonexistence for $\mu>\LN$ in the class of solutions of Theorem \ref{thm:main1} (which do not belong to $H^1_{\rm loc}(\Omega)$ in general) remain open.

 The uniqueness argument of Theorem \ref{thm:uniqueness} follows the ideas of \cite{boccardo2014some, canino2004variational} for the local case; the nonexistence argument of Theorem \ref{thm:nonexistence} combines the local Picone identity \cite{allegretto1998picone} with its nonlocal counterpart \cite{brasco2014convexity, frank2008nonlinear}; the lower bound of Theorem \ref{thm:origin} extends to $\cL$ a classical comparison argument for the local operator \cite{brezis2005semilinear, dupaigne2002nonlinear}, the main new ingredient being the control of the nonlocal part of the explicit subsolution (Lemma \ref{lem:subsol}), for which we use the homogeneity of $(-\Delta)^s|x|^{-\beta}$ (whose exact value is computed in \cite{frank2008hardy}); the critical case relies on the improved Hardy inequalities of \cite{brezis1997blow, vazquez2000hardy}. What we believe is new is the treatment of the combined problem \eqref{problem1}: the passage to the limit in the nonlocal term for solutions of infinite energy (Lemma \ref{lem:welldef}), the explicit thresholds $\mu_{m,\gamma}$ and their optimality for $\delta_m\ge1$ (Theorem \ref{thm:origin}(c)), and the results of Theorems \ref{thm:summability}--\ref{thm:critical} for the mixed operator.

The paper is organised as follows. In Section \ref{sec:prelim} we collect the functional setting and several auxiliary results. Section \ref{sec:approx} is devoted to the approximating problems, Section \ref{sec:estimates} to a priori estimates, and Section \ref{sec:proofs} to the proofs of Theorems \ref{thm:main1}--\ref{thm:main3}. In Section \ref{sec:further} we prove Theorems \ref{thm:summability}--\ref{thm:nonexistence}; the main new tools there are a test-function extension argument for the energy formulation and Picone-type inequalities for the mixed operator. Section \ref{sec:origin} contains the proofs of Theorems \ref{thm:origin} and \ref{thm:critical}, based on an explicit subsolution of the form $|x|^{-\alpha_\mu}+|x|^{-\alpha_\mu+\tau}$ and on a comparison argument.

\section{Preliminaries}\label{sec:prelim}
In this section, we present the basic definitions and results concerning the function spaces, inequalities, and elementary estimates used throughout the paper. To make the paper self-contained, we also provide proofs of several results in the appendix, although some of them have already appeared in earlier works. We begin by setting up the functional space for our study. 

For $1\le r\le\infty$ we denote by $\|\cdot\|_r$ the norm of $L^r(\Omega)$. Every function defined in $\Omega$ is tacitly extended by zero to $\RN\setminus\Omega$. We let $H^1_0(\Omega)$ be the closure of $C^\infty_c(\Omega)$ in $H^1(\RN)$, endowed with the norm $\|\nabla u\|_2$, and for $1\le q<\infty$ we define $W^{1,q}_0(\Omega)$ analogously. For $u\in H^1(\RN)$ the Gagliardo seminorm of order $s$ is
\[
[u]_s:=\Big(\iint_{\R^{2N}}\frac{|u(x)-u(y)|^2}{|x-y|^{N+2s}}\dd x\dd y\Big)^{1/2},
\]
and we set
\[
\rho(u):=\cE(u,u)^{1/2}=\Big(\|\nabla u\|_{L^2(\RN)}^2+\frac{\Cs}{2}[u]_s^2\Big)^{1/2}.
\]
By \cite[Proposition 2.2]{dinezza2012hitchhiker}, $[u]_s\le C\|u\|_{H^1(\RN)}$ for every $u\in H^1(\RN)$; combining this with the Poincar\'e inequality in $H^1_0(\Omega)$ we obtain a constant $C_*=C_*(N,s,\Omega)\ge1$ such that
\begin{equation}\label{eq:normequiv}
\|\nabla u\|_2\le\rho(u)\le C_*\|\nabla u\|_2\qquad\text{for all }u\in H^1_0(\Omega).
\end{equation}
Thus $\rho$ is a norm on $H^1_0(\Omega)$ equivalent to $\|\nabla\cdot\|_2$; in particular, $H^1_0(\Omega)$ coincides with the completion of $C^\infty_c(\Omega)$ with respect to $\rho$, which is the natural energy space for $\cL$ (see \cite{biagi2022mixed, biagi2025mixed}). By the Cauchy--Schwarz inequality, $|\cE(u,\psi)|\le\rho(u)\rho(\psi)$ for $u,\psi\in H^1_0(\Omega)$.

We denote by $S_N$ the constant of the Sobolev inequality and by $C_P$ the Poincar\'e constant of $\Omega$:
\begin{equation}\label{eq:sobolev}
\|u\|_{2^*}^2\le S_N\|\nabla u\|_2^2,\qquad \|u\|_2\le C_P\|\nabla u\|_2,\qquad u\in H^1_0(\Omega).
\end{equation}
For $k>0$ we use the truncations $T_k(t):=\min\{t,k\}$ and $G_k(t):=(t-k)^+$, $t\ge0$. Moreover, for a function $z$ we write $z^+:=\max\{z,0\}$ and $z^-:=\max\{-z,0\}$, so that $z=z^+-z^-$.

The next lemma collects the elementary properties of the exponents \eqref{eq:msharp}--\eqref{eq:threshold}.

\begin{lemma}\label{lem:exponents}
Let $(m,\gamma)\in\cS$ with $1\le m<N/2$. Then:
\begin{enumerate}[label=(\roman*)]
\item $\displaystyle 1+\delta_m=\frac{m(N-2)(1+\gamma)}{N-2m}$;
\item $m\mapsto\delta_m$ is increasing on $[1,N/2)$, $\delta_m>0$, and $\delta_m\ge\gamma$ with equality if and only if $m=1$. If $0\le\gamma<1$, then $m^\#<N/2$ and $\delta_m<1$ if and only if $m<m^\#$; if $\gamma\ge1$, then $\delta_m\ge1$;
\item if $m>1$, then $\displaystyle(\delta_m-\gamma)m'=\frac{2^*}{2}(1+\delta_m)=\sigma_m$;
\item if $0\le\gamma<1$ and $m<m^\#$, then $\displaystyle\frac{(1-\delta_m)q_m}{2-q_m}=\sigma_m$, $\ \dfrac{N}{N-1}<q_m<2$, and $q_m^*:=\dfrac{Nq_m}{N-q_m}=\sigma_m$;
\item the two expressions in \eqref{eq:threshold} coincide, and $0<\mu_{m,\gamma}\le\LN$, with equality if and only if $\delta_m=1$.
\end{enumerate}
\end{lemma}

\begin{proof}
(i) We have
\[
1+\delta_m=\frac{N(m-1)+\gamma m(N-2)+N-2m}{N-2m}=\frac{m(N-2)+\gamma m(N-2)}{N-2m}.
\]

(ii) The map $m\mapsto\delta_m$ is increasing on $[1,N/2)$, since both its numerator is increasing and its (positive) denominator is decreasing in $m$, and $\delta_1=\gamma$. If $\gamma>0$ this gives $\delta_m\ge\gamma>0$; if $\gamma=0$, then $m>1$ because $(1,0)\notin\cS$, and $\delta_m=N(m-1)/(N-2m)>0$. If $\gamma\ge1$, then $\delta_m\ge\delta_1=\gamma\ge1$. If $0\le\gamma<1$, then $m^\#\le 2N/(N+2)<N/2$ because $N\ge3$, and
\[
\delta_m<1\iff N(m-1)+\gamma m(N-2)<N-2m\iff m\big(N+2+\gamma(N-2)\big)<2N ,
\]
that is, $\delta_m<1$ if and only if $m<m^\#$.

(iii) The identity $(\delta-\gamma)\frac{m}{m-1}=\frac{N}{N-2}(1+\delta)$ is equivalent to $\delta\,(N-2m)=N(m-1)+\gamma m(N-2)$, which is the definition of $\delta_m$. The equality $\frac{N}{N-2}(1+\delta_m)=\sigma_m$ follows from (i).

(iv) By (i),
\[
1-\delta_m+\sigma_m=2+\frac{Nm(1+\gamma)-m(N-2)(1+\gamma)}{N-2m}=\frac{2\big(N-m(1-\gamma)\big)}{N-2m},
\]
hence $q_m=2\sigma_m/(1-\delta_m+\sigma_m)$, which is equivalent to $(1-\delta_m)q_m=\sigma_m(2-q_m)$. Since $\delta_m<1$ we get $q_m<2$. The map $m\mapsto q_m$ is increasing; if $\gamma>0$, then
\[
q_m\ge q_1=\frac{N(1+\gamma)}{N-1+\gamma}>\frac{N}{N-1}\iff (1+\gamma)(N-1)>N-1+\gamma\iff\gamma(N-2)>0,
\]
while if $\gamma=0$ then $m>1$ and $q_m=Nm/(N-m)>N/(N-1)$. Finally $N-q_m=N(N-2m)/(N-m(1-\gamma))$, so $Nq_m/(N-q_m)=\sigma_m$.

(v) Using (i) and the definition of $\delta_m$,
\[
\frac{4\delta_m}{(1+\delta_m)^2}\LN=\frac{N(m-1)+\gamma m(N-2)}{N-2m}\cdot\frac{(N-2m)^2}{m^2(N-2)^2(1+\gamma)^2}\cdot(N-2)^2 ,
\]
which is the second expression in \eqref{eq:threshold}. Since $\delta_m>0$ and $(1+\delta_m)^2-4\delta_m=(1-\delta_m)^2\ge0$, we have $0<4\delta_m/(1+\delta_m)^2\le1$, with equality if and only if $\delta_m=1$.
\end{proof}

 We now recall the following Hardy inequalities; see \cite{garcia1998hardy, peral2021elliptic}.

\begin{lemma}\label{lem:hardy}
\begin{enumerate}[label=(\alph*)]
\item For every $u\in H^1_0(\Omega)$,
\begin{equation}\label{eq:hardy}
\LN\int_\Omega\frac{u^2}{|x|^2}\dd x\le\int_\Omega|\nabla u|^2\dd x\le\rho(u)^2 .
\end{equation}
\item Let $1\le q<N$. For every $v\in W^{1,q}(\RN)$ with compact support,
\begin{equation}\label{eq:hardyq}
\int_{\RN}\frac{|v|^q}{|x|^q}\dd x\le\Big(\frac{q}{N-q}\Big)^q\int_{\RN}|\nabla v|^q\dd x .
\end{equation}
\end{enumerate}
\end{lemma}


\begin{remark}\label{rem:hardy-optimal}
The constant $\LN$ in \eqref{eq:hardy} is optimal also with respect to $\rho$, i.e.
\[
\inf\Big\{\rho(u)^2\Big/\!\int_\Omega\frac{u^2}{|x|^2}\dd x:\ u\in H^1_0(\Omega)\setminus\{0\}\Big\}=\LN
\]
(see \cite{biagi2025mixed} for a thorough discussion). Indeed, given $u\in C^\infty_c(\RN)\setminus\{0\}$, the function $u_\lambda(x):=u(\lambda x)$ is supported in $\Omega$ for $\lambda$ large, because $0\in\Omega$, and
\[
\frac{\rho(u_\lambda)^2}{\int u_\lambda^2|x|^{-2}\dd x}=\frac{\|\nabla u\|^2_{L^2(\RN)}+\frac{\Cs}{2}\lambda^{2s-2}[u]_s^2}{\int u^2|x|^{-2}\dd x}\longrightarrow\frac{\|\nabla u\|^2_{L^2(\RN)}}{\int u^2|x|^{-2}\dd x}\qquad(\lambda\to\infty),
\]
and the infimum of the right-hand side over $u$ is $\LN$.
\end{remark}

The following algebraic inequalities will be used repeatedly throughout the paper (see e.g. \cite{arora2023combined}). 
\begin{lemma}\label{lem:alg}
\begin{enumerate}[label=(\roman*)]
\item For all $a,b\ge0$ and $\theta>0$,
\begin{equation}\label{eq:alg}
\frac{4\theta}{(\theta+1)^2}\Big(a^{\frac{\theta+1}{2}}-b^{\frac{\theta+1}{2}}\Big)^2\le(a-b)(a^\theta-b^\theta).
\end{equation}
\item If $\Phi:\R\to\R$ is nondecreasing, then $(a-b)(\Phi(a)-\Phi(b))\ge0$ for all $a,b\in\R$. If, moreover, $\Phi$ is $1$-Lipschitz, then $(a-b)(\Phi(a)-\Phi(b))\ge(\Phi(a)-\Phi(b))^2$.
\item For every $z\in H^1_0(\Omega)$ we have $z^\pm\in H^1_0(\Omega)$ and
\[
\cE(z,z^+)\ge\rho(z^+)^2,\qquad \cE(z,-z^-)\ge\rho(z^-)^2 .
\]
\end{enumerate}
\end{lemma}

In the following lemma, a weak solution of \[ \cL w=f \quad \text{in }\Omega, \qquad w=0 \quad \text{in }\RN\setminus\Omega, \] with $f\in L^{(2^*)'}(\Omega)$, is a function $w\in H^1_0(\Omega)$ such that \[ \cE(w,\psi)=\int_\Omega f\psi\dd x \] for every $\psi\in H^1_0(\Omega)$. The proof is based on the classical Stampacchia argument \cite{stampacchia1965dirichlet}, combined with the nonnegativity of the nonlocal part of the form when tested against nondecreasing functions of the solution. Related arguments can be found in \cite{biagi2022mixed,leonori2015basic}. We provide the details in the appendix for completeness.

\begin{lemma}\label{lem:stampacchia} Let $f\in L^p(\Omega)$ with $p\ge (2^*)'$ and $f\ge 0$, and let $w\in H^1_0(\Omega)$, $w\ge 0$, be a weak solution of 

$$ \cL w=f \quad \text{in }\Omega, \qquad w=0 \quad \text{in }\RN\setminus\Omega. $$

Then the following assertions hold. 

\begin{enumerate}[label=(\alph*)] \item If $p>N/2$, then $w\in L^\infty(\Omega)$ and $$ \|w\|_\infty\le C(N,p,|\Omega|)\|f\|_p. $$ 
\item If $(2^*)'\le p<N/2$, then $w\in L^{p^{**}}(\Omega)$, where $$ p^{**}:=\frac{Np}{N-2p}, $$ and $$ \|w\|_{p^{**}}\le C(N,p)\|f\|_p. $$ \end{enumerate} \end{lemma}

 We say that $w\in H^1_0(\Omega)$ is a weak supersolution of $$ \cL w=0 \quad \text{in }\Omega $$ if $\cE(w,\varphi)\ge 0 $ for every nonnegative $\varphi\in C^\infty_c(\Omega)$. By density and the continuity of the form $\cE$, the same inequality holds for every nonnegative $\varphi\in H^1_0(\Omega)$ with compact support in $\Omega$. 

\begin{lemma}\label{lem:strongmin} Let $w\in H^1_0(\Omega)$ be a weak supersolution of $$ \cL w=0 \quad \text{in }\Omega $$ such that $w\ge 0$ a.e. in $\RN$ and $w\not\equiv 0$ in $\Omega$. Then, for every $\omega\Subset\Omega$, there exists $c_\omega>0$ such that $$ w\ge c_\omega \quad\text{a.e. in }\omega. $$ \end{lemma} 

 The following lemma shows that Definition~\ref{def:weak} is meaningful and provides the convergence result for the nonlocal term needed when passing to the limit. 

\begin{lemma}\label{lem:welldef}
Let $q\in(\frac{N}{N-1},N)$ and let $u\in L^1(\Omega)\cap W^{1,q}_{\rm loc}(\Omega)$, with $u=0$ in $\RN\setminus\Omega$. Then, for every $\psi\in C^\infty_c(\Omega)$:
\begin{enumerate}[label=(\roman*)]
\item $\nabla u\cdot\nabla\psi\in L^1(\Omega)$ and $u\psi|x|^{-2}\in L^1(\Omega)$;
\item the function $(x,y)\mapsto(u(x)-u(y))(\psi(x)-\psi(y))|x-y|^{-N-2s}$ belongs to $L^1(\R^{2N})$, so that $\cE(u,\psi)$ is well defined;
\item if $(u_n)$ is a sequence of functions with the same properties as $u$ such that $u_n\to u$ in $L^1(\Omega)$ and $u_n\rightharpoonup u$ weakly in $W^{1,q}(U)$ for every open set $U\Subset\Omega$, then $\cE(u_n,\psi)\to\cE(u,\psi)$.
\end{enumerate}
\end{lemma}

\begin{proof}  Fix $\psi\in C^\infty_c(\Omega)$ and set $$ K:=\supp\psi, \qquad d:=\frac13 {\rm dist}(K,\partial\Omega)>0. $$ Define $ U':=\{x\in\RN:{\rm dist}(x,K)<d\} $ and $$ U:=\{x\in\RN:{\rm dist}(x,K)<2d\}. $$ Then $ K\Subset U'\Subset U\Subset\Omega. $ Choose $\varphi\in C^\infty_c(U)$ such that $0\le\varphi\le 1$ and $\varphi\equiv 1$ in $U'$, and write $$ u=u\varphi+u(1-\varphi)=:u_1+u_2. $$ Since $u\in W^{1,q}(U)$, the function $u_1$, extended by zero outside $U$, belongs to $W^{1,q}(\RN)$ and has compact support in $U$. Moreover, $$ \|u_1\|_{W^{1,q}(\RN)} \le C_\varphi\|u\|_{W^{1,q}(U)}. $$

(i)  Since $\supp\nabla\psi\subset K$ and $u\in W^{1,q}(U)$, H\"older's inequality gives $$ \int_\Omega|\nabla u\cdot\nabla\psi|\dd x \le \|\nabla u\|_{L^q(K)} \|\nabla\psi\|_{L^{q'}(K)} <\infty. $$ 

Thus $\nabla u\cdot\nabla\psi\in L^1(\Omega)$. Since $\varphi\equiv 1$ on $K$, we have $u\psi=u_1\psi$. By the Hardy inequality \eqref{eq:hardyq}, 
$$ \left\|\frac{u_1}{|x|}\right\|_{L^q(\RN)} \le C\|\nabla u_1\|_{L^q(\RN)}. $$ 

Furthermore, the assumption $q>N/(N-1)$ implies $q'<N$, and hence $ \frac{\psi}{|x|}\in L^{q'}(\RN). $ Therefore, H\"older's inequality yields $$ \int_\Omega\frac{|u\psi|}{|x|^2}\dd x = \int_{\RN} \frac{|u_1|}{|x|} \frac{|\psi|}{|x|}\dd x \le \left\|\frac{u_1}{|x|}\right\|_q \left\|\frac{\psi}{|x|}\right\|_{q'} <\infty. $$

\medskip

(ii) Choose $t\in(0,1)$ such that $ 0<2s-t<1. $ Such a choice is possible because $ \max\{0,2s-1\}<\min\{1,2s\}. $ Using $ |x-y|^{N+2s} = |x-y|^{\frac Nq+t} |x-y|^{\frac N{q'}+(2s-t)} $ and applying H\"older's inequality in $\R^{2N}$, we obtain \begin{equation}\label{eq:holderfrac} \begin{aligned} &\iint_{\R^{2N}} \frac{|u_1(x)-u_1(y)||\psi(x)-\psi(y)|} {|x-y|^{N+2s}}\dd x\dd y \le [u_1]_{W^{t,q}(\RN)} [\psi]_{W^{2s-t,q'}(\RN)}, \end{aligned} \end{equation} where, for $\sigma\in(0,1)$ and $1\le r<\infty$, $$ [v]_{W^{\sigma,r}(\RN)} := \left( \iint_{\R^{2N}} \frac{|v(x)-v(y)|^r} {|x-y|^{N+\sigma r}}\dd x\dd y \right)^{1/r}. $$ 

Both seminorms on the right-hand side of \eqref{eq:holderfrac} are finite. Indeed, $ [\psi]_{W^{2s-t,q'}(\RN)}<\infty $ because $\psi\in C^\infty_c(\RN)$ and $0<2s-t<1$, while \cite[Proposition 2.2]{dinezza2012hitchhiker} gives $ [u_1]_{W^{t,q}(\RN)} \le C\|u_1\|_{W^{1,q}(\RN)}. $ It remains to consider $u_2$. Since $u_2=0$ in $U'$ and $\psi=0$ outside $K$, the product $ (u_2(x)-u_2(y))(\psi(x)-\psi(y)) $ can be nonzero only if one of the two points lies in $K$ and the other lies outside $U'$. If, for example, $x\in K$ and $y\in\RN\setminus U'$, then $|x-y|\ge d$ and $$ (u_2(x)-u_2(y))(\psi(x)-\psi(y)) = -u_2(y)\psi(x). $$ 

By symmetry, \begin{equation}\label{eq:farpart} \begin{aligned} &\frac{\Cs}{2} \iint_{\R^{2N}} \frac{(u_2(x)-u_2(y))(\psi(x)-\psi(y))} {|x-y|^{N+2s}}\dd x\dd y = -\Cs\int_K\psi(x) \int_{\RN\setminus U'} \frac{u_2(y)} {|x-y|^{N+2s}}\dd y\dd x. \end{aligned} \end{equation} Moreover, the absolute value of the integrand on the right-hand side is bounded by $$ d^{-N-2s}\|\psi\|_\infty \chi_K(x)|u(y)|. $$ This function is integrable on $\R^{2N}$ because $u=0$ outside $\Omega$ and $u\in L^1(\Omega)$. Combining this observation with \eqref{eq:holderfrac} proves \textup{(ii)}.

\medskip

(iii) Set $ u_{n,1}:=u_n\varphi, \, u_{n,2}:=u_n(1-\varphi). $ Since $ u_n\rightharpoonup u$ weakly in $W^{1,q}(U),$ multiplication by $\varphi$ gives $ u_{n,1}\rightharpoonup u_1$ weakly in $W^{1,q}(\RN).$  In particular, $$ \int_\Omega\nabla u_n\cdot\nabla\psi\dd x \to \int_\Omega\nabla u\cdot\nabla\psi\dd x. $$ 

We next consider the Hardy-potential term. For $v\in W^{1,q}(\RN)$, the Hardy inequality and H\"older's inequality give $$ \begin{aligned} \left| \int_{\RN}\frac{v\psi}{|x|^2}\dd x \right| \le \left\|\frac{v}{|x|}\right\|_q \left\|\frac{\psi}{|x|}\right\|_{q'} \le C \|v\|_{W^{1,q}(\RN)} \left\|\frac{\psi}{|x|}\right\|_{q'}. \end{aligned} $$ 

Thus $ v\mapsto \int_{\RN}\frac{v\psi}{|x|^2}\dd x $ is a bounded linear functional on $W^{1,q}(\RN)$. Since $u_n\psi=u_{n,1}\psi$, it follows that $$ \int_\Omega\frac{u_n\psi}{|x|^2}\dd x \to \int_\Omega\frac{u\psi}{|x|^2}\dd x. $$ By \eqref{eq:holderfrac}, the map $$ v\longmapsto \iint_{\R^{2N}} \frac{(v(x)-v(y))(\psi(x)-\psi(y))} {|x-y|^{N+2s}}\dd x\dd y $$ is a bounded linear functional on $W^{1,q}(\RN)$. It is therefore weakly continuous. Hence the nonlocal part of $\cE(u_{n,1},\psi)$ converges to the nonlocal part of $\cE(u_1,\psi)$. Finally, applying \eqref{eq:farpart} to $u_{n,2}-u_2$, we obtain $$ \begin{aligned} &\left| \frac{\Cs}{2} \iint_{\R^{2N}} \frac{ \bigl((u_{n,2}-u_2)(x)-(u_{n,2}-u_2)(y)\bigr) (\psi(x)-\psi(y)) } {|x-y|^{N+2s}}\dd x\dd y \right| \\ &\qquad\le \Cs d^{-N-2s}\|\psi\|_{L^1(\Omega)} \|u_n-u\|_{L^1(\Omega)}. \end{aligned} $$ The right-hand side tends to zero. Combining the convergence of the local, Hardy-potential, and nonlocal terms yields $ \cE(u_n,\psi)\to\cE(u,\psi). $
\end{proof}

\begin{remark}\label{rem:welldef}
Suppose that $u$ satisfies conditions \textup{(a)} and \textup{(b)} of Definition~\ref{def:weak}. Let $\psi\in C^\infty_c(\Omega)$, $g\in L^1(\Omega)$, and $K:=\supp\psi$. If $ u\ge c_K>0$ a.e. in $K,$ then $ |g\psi u^{-\gamma}| \le c_K^{-\gamma}\|\psi\|_\infty |g|$ a.e. in $\Omega. $ Consequently, $$ \int_\Omega|g\psi u^{-\gamma}|\dd x \le c_K^{-\gamma}\|\psi\|_\infty\|g\|_1<\infty, $$ so the right-hand side of \eqref{eq:weakform} is well defined.
\end{remark}

\section{The approximating problems}\label{sec:approx}

The approximation scheme of this section is the one of Boccardo and Orsina \cite{boccardo2010semilinear} (see also \cite{arora2023combined, garain2023class} for the mixed operator without Hardy term), with the additional regularization $|x|^{-2}\mapsto(|x|^2+1/n)^{-1}$ of the potential, as in \cite{boccardo2006existence}. Throughout this section $g\in L^1(\Omega)$ is nonnegative and not identically zero, $\gamma\ge0$, $0<\mu<\LN$, and for $n\in\N$ we set
\[
g_n:=\min\{g,n\},\qquad V_n(x):=\frac{1}{|x|^2+1/n}.
\]
Note that $0\le g_n\le g_{n+1}\le g$, $g_n\to g$ a.e. in $\Omega$ (and in $L^m(\Omega)$ if $g\in L^m(\Omega)$, $m<\infty$), and
\begin{equation}\label{eq:Vn}
0<V_n\le V_{n+1}\le\frac{1}{|x|^2},\qquad V_n\le n .
\end{equation}
A weak solution of \eqref{problem11} is a nonnegative function $u_n\in H^1_0(\Omega)$ such that
\begin{equation}\label{weak11}
\cE(u_n,\psi)-\mu\int_\Omega V_n u_n\psi\dd x=\int_\Omega\frac{g_n\,\psi}{(u_n+1/n)^\gamma}\dd x\qquad\text{for all }\psi\in H^1_0(\Omega).
\end{equation}

\subsection{A linear problem}
For $w,\psi\in H^1_0(\Omega)$ set
\[
\cB_n(w,\psi):=\cE(w,\psi)-\mu\int_\Omega V_n w\psi\dd x .
\]

\begin{lemma}\label{lem:linear}
Let $n\in\N$ and $h\in L^2(\Omega)$. There exists a unique $w\in H^1_0(\Omega)$ such that
\begin{equation}\label{weak_S_n}
\cB_n(w,\psi)=\int_\Omega h\psi\dd x\qquad\text{for all }\psi\in H^1_0(\Omega),
\end{equation}
and it satisfies $\big(1-\frac{\mu}{\LN}\big)\|\nabla w\|_2\le C_P\|h\|_2$. If, moreover, $h\ge0$ a.e. in $\Omega$, then $w\ge0$ a.e. in $\RN$.
\end{lemma}

\begin{proof}
By \eqref{eq:Vn}, the Cauchy--Schwarz inequality and \eqref{eq:hardy},
\begin{equation}\label{eq:hardyterm-bound}
\mu\Big|\int_\Omega V_nw\psi\dd x\Big|\le\mu\Big\|\frac{w}{|x|}\Big\|_2\Big\|\frac{\psi}{|x|}\Big\|_2\le\frac{\mu}{\LN}\|\nabla w\|_2\|\nabla\psi\|_2 .
\end{equation}
Together with $|\cE(w,\psi)|\le\rho(w)\rho(\psi)$ and \eqref{eq:normequiv}, this shows that $\cB_n$ is a bounded bilinear form on $H^1_0(\Omega)$. Moreover, for every $u\in H^1_0(\Omega)$, again by \eqref{eq:Vn},  the Hardy inequality \eqref{eq:hardy} and \eqref{eq:normequiv},
\begin{equation}\label{coercive}
\cB_n(u,u)\ge\rho(u)^2-\mu\int_\Omega\frac{u^2}{|x|^2}\dd x\ge\rho(u)^2-\frac{\mu}{\LN}\|\nabla u\|_2^2\ge\Big(1-\frac{\mu}{\LN}\Big)\rho(u)^2\ge\Big(1-\frac{\mu}{\LN}\Big)\|\nabla u\|_2^2,
\end{equation}
so $\cB_n$ is coercive because $\mu<\LN$. The linear functional $\psi\mapsto\int_\Omega h\psi\dd x$ is bounded on $H^1_0(\Omega)$, since $$\left|\int_\Omega h\psi\right| \le\|h\|_2\|\psi\|_2\le C_P\|h\|_2\|\nabla\psi\|_2.$$ By the Lax--Milgram theorem \cite[Corollary 5.8]{Brezis} there is a unique $w\in H^1_0(\Omega)$ satisfying \eqref{weak_S_n}. Testing \eqref{weak_S_n} with $\psi=w$ and using \eqref{coercive} we get $$(1-\mu/\LN)\|\nabla w\|_2^2\le C_P\|h\|_2\|\nabla w\|_2.$$

Assume now $h\ge0$, and test \eqref{weak_S_n} with $\psi=-w^-\in H^1_0(\Omega)$, where $w^-=\max\{-w,0\}$. Since $V_nw(-w^-)=V_n(w^-)^2$, Lemma \ref{lem:alg}(iii) and \eqref{coercive} give
\[
\Big(1-\frac{\mu}{\LN}\Big)\rho(w^-)^2\le\cB_n(w^-,w^-)\le\cB_n(w,-w^-)=-\int_\Omega h\,w^-\dd x\le0 .
\]
Hence $w^-=0$, i.e. $w\ge0$ a.e.
\end{proof}

\subsection{\texorpdfstring{Existence, uniqueness and boundedness of $u_n$}{Existence, uniqueness and boundedness}}

\begin{lemma}\label{lem:Pn}
For every $n\in\N$, problem \eqref{problem11} has a unique weak solution $u_n$. Moreover $u_n\in L^\infty(\Omega)$.
\end{lemma}

\begin{proof} We divide the proof into three steps.

\emph{Existence.} For $\xi\in L^2(\Omega)$, set $$ h_\xi:=\frac{g_n}{(\xi^++1/n)^\gamma}. $$ Since $0\le g_n\le n$, we have $ 0\le h_\xi\le n^{1+\gamma}. $ Let $T(\xi):=w$, where $w\in H^1_0(\Omega)$ is the unique solution of \eqref{weak_S_n} with $h=h_\xi$, whose existence is provided by Lemma~\ref{lem:linear}. By the estimate in Lemma~\ref{lem:linear}, \begin{equation}\label{radius} \|\nabla T(\xi)\|_2 \le R_n := \frac{C_Pn^{1+\gamma}|\Omega|^{1/2}} {1-\mu/\LN} \end{equation} for every $\xi\in L^2(\Omega)$. We regard $T$ as a map from $L^2(\Omega)$ into itself and consider the closed convex set $ \mathcal K := \left\{ \xi\in L^2(\Omega): \|\xi\|_2\le C_PR_n \right\}. $ By \eqref{radius} and the Poincar\'e inequality, $$ \|T(\xi)\|_2 \le C_P\|\nabla T(\xi)\|_2 \le C_PR_n. $$ Thus, $ T\bigl(L^2(\Omega)\bigr)\subset\mathcal K, $ and, in particular, $T(\mathcal K)\subset\mathcal K$. Moreover, \eqref{radius} shows that $T(\mathcal K)$ is bounded in $H^1_0(\Omega)$. Hence $T(\mathcal K)$ is relatively compact in $L^2(\Omega)$ by the Rellich--Kondrachov theorem. 

It remains to prove that $T$ is continuous. Let $ \xi_k\to\xi$ in $L^2(\Omega), $ and set $ w_k:=T(\xi_k),$ \, $w:=T(\xi). $ Subtracting the weak formulations satisfied by $w_k$ and $w$, testing the resulting identity with $w_k-w$, and using \eqref{coercive} and the Poincar\'e inequality, we obtain \begin{equation}\label{subtractw} \begin{aligned} \left(1-\frac{\mu}{\LN}\right) \|\nabla(w_k-w)\|_2^2 &\le \int_\Omega (h_{\xi_k}-h_\xi)(w_k-w)\dd x \\ &\le C_P\|h_{\xi_k}-h_\xi\|_2 \|\nabla(w_k-w)\|_2. \end{aligned} \end{equation}

We claim that $ h_{\xi_k}\to h_\xi$ in $L^2(\Omega). $ Indeed, every subsequence of $(\xi_k)$ has a further subsequence converging to $\xi$ a.e. in $\Omega$. Along such a subsequence, the continuity of $ t\longmapsto (t^++1/n)^{-\gamma} $ gives $ h_{\xi_k}\to h_\xi$ a.e. in  $\Omega.$

Moreover, $ |h_{\xi_k}-h_\xi| \le n^{1+\gamma}. $ Since $|\Omega|<\infty$, the dominated convergence theorem gives convergence in $L^2(\Omega)$ along the extracted subsequence. As this argument applies to every subsequence, the full sequence satisfies $ h_{\xi_k}\to h_\xi $ in $L^2(\Omega). $ It now follows from \eqref{subtractw} that $ w_k\to w$ in $H^1_0(\Omega), $ and hence in $L^2(\Omega)$. Therefore, $T:\mathcal K\to\mathcal K$ is continuous and compact. By the Schauder fixed point theorem in the form of \cite[Corollary 11.2]{gilbarg2001elliptic}, there exists $u_n\in\mathcal K$ such that $ T(u_n)=u_n. $ By Lemma~\ref{lem:linear}, $u_n\ge0$ a.e. in $\Omega$. Hence $u_n^+=u_n$, and the fixed-point identity shows that $u_n$ is a weak solution of \eqref{problem11}.

\emph{Uniqueness.} Let $u$ and $\widetilde u$ be two weak solutions of \eqref{problem11}, and set $ z:=u-\widetilde u. $ Subtracting the corresponding identities \eqref{weak11} and testing with $z^+$, we obtain $$ \cB_n(z,z^+) = \int_\Omega g_n \left( \frac{1}{(u+1/n)^\gamma} - \frac{1}{(\widetilde u+1/n)^\gamma} \right) z^+\dd x. $$ The function $ t\longmapsto \frac{1}{(t+1/n)^\gamma} $ is decreasing on $[0,\infty)$. Therefore, the integrand on the right-hand side is nonpositive on the set $\{z>0\}$, and hence $ \cB_n(z,z^+)\le0. $ On the other hand, we have $ V_nzz^+=V_n(z^+)^2$ a.e. in $\Omega. $ Together with Lemma~\ref{lem:alg}\textup{(iii)} and \eqref{coercive}, this gives $$ \cB_n(z,z^+) \ge \cB_n(z^+,z^+) \ge \left(1-\frac{\mu}{\LN}\right)\rho(z^+)^2. $$

Since $\mu<\LN$, it follows that $ \rho(z^+)^2=0, $ and therefore $z^+=0$. Thus, $ u\le\widetilde u \quad\text{a.e. in }\Omega. $ Exchanging the roles of $u$ and $\widetilde u$, we obtain the reverse inequality. Consequently, $ u=\widetilde u \quad\text{a.e. in }\Omega. $

\emph{Boundedness.} The function $u_n$ is a weak solution, in the sense of Lemma~\ref{lem:stampacchia}, of $ \cL u_n=f_n \quad\text{in }\Omega, $ where $$ f_n := \mu V_nu_n + \frac{g_n}{(u_n+1/n)^\gamma}. $$ 

By \eqref{eq:Vn}, the nonnegativity of $u_n$, and the bounds $0\le g_n\le n$, we have $ 0\le f_n \le \mu n u_n+n^{1+\gamma}. $ Since $|\Omega|<\infty$, it follows that $ u_n\in L^p(\Omega) \Longrightarrow f_n\in L^p(\Omega) $ for every $p\ge1$. By the Sobolev inequality, $ u_n\in L^{2^*}(\Omega). $ Set $ p_0:=2^*. $ Then $f_n\in L^{p_0}(\Omega)$ and $p_0>(2^*)'$. If $ p_0>\frac N2, $ Lemma~\ref{lem:stampacchia}\textup{(a)} immediately gives $ u_n\in L^\infty(\Omega). $ Suppose now that $p_0\le N/2$. As long as $p_j<N/2$, define $p_{j+1}$ by $$ \frac{1}{p_{j+1}} := \frac{1}{p_j}-\frac2N. $$ Equivalently, $ p_{j+1}=p_j^{**} = \frac{Np_j}{N-2p_j}. $ If $u_n\in L^{p_j}(\Omega)$, then $f_n\in L^{p_j}(\Omega)$, and Lemma~\ref{lem:stampacchia}\textup{(b)} gives $ u_n\in L^{p_{j+1}}(\Omega). $ Since $$ \frac1{p_j} = \frac1{p_0}-\frac{2j}{N}, $$ after finitely many steps we reach an index $j_0$ such that $ p_{j_0}\ge\frac N2. $ If $p_{j_0}>N/2$, then $f_n\in L^{p_{j_0}}(\Omega)$ and Lemma~\ref{lem:stampacchia}\textup{(a)} yields $ u_n\in L^\infty(\Omega). $ It remains to consider the case $ p_{j_0}=\frac N2. $ Since $\Omega$ is bounded and $ u_n,f_n\in L^{N/2}(\Omega), $ we also have $ u_n,f_n\in L^p(\Omega) $ for every $1\le p<N/2$. 

Choose $$ p\in \left( \max\left\{\frac N4,(2^*)'\right\}, \frac N2 \right). $$ This interval is nonempty. Lemma~\ref{lem:stampacchia}\textup{(b)} then gives $ u_n\in L^{p^{**}}(\Omega), \,\, p^{**}=\frac{Np}{N-2p}. $ Since $p>N/4$, we have $ p^{**}>\frac N2. $ Consequently, $ f_n\in L^{p^{**}}(\Omega), $ and a final application of Lemma~\ref{lem:stampacchia}\textup{(a)} gives $ u_n\in L^\infty(\Omega). $ This completes the proof.
\end{proof}

\subsection{Monotonicity and uniform positivity}

\begin{lemma}\label{lem:monotone} The sequence $(u_n)$ is nondecreasing, that is, $u_n\le u_{n+1}$ a.e. in $\Omega$ for every $n\in\N$. Moreover, for every $\omega\Subset\Omega$, there exists $c_\omega>0$, independent of $n$, such that 

\begin{equation}\label{eq:unifpos} u_n\ge c_\omega \qquad\text{a.e. in }\omega,\quad\text{for every }n\in\N. \end{equation} \end{lemma} 

\begin{proof} Fix $n\in\N$ and set $z:=u_n-u_{n+1}\in H^1_0(\Omega)$. Subtracting the identities \eqref{weak11} satisfied by $u_n$ and $u_{n+1}$ and testing with $z^+$, we obtain \begin{equation}\label{delta1} \begin{aligned} &\cE(z,z^+) -\mu\int_\Omega \bigl(V_nu_n-V_{n+1}u_{n+1}\bigr)z^+\dd x = \int_\Omega \left( \frac{g_n}{(u_n+1/n)^\gamma} -\frac{g_{n+1}}{(u_{n+1}+1/(n+1))^\gamma} \right)z^+\dd x. \end{aligned} \end{equation}

We estimate the terms in \eqref{delta1}. First, Lemma~\ref{lem:alg}\textup{(iii)} gives $ \cE(z,z^+)\ge\rho(z^+)^2. $ Second, since $u_n\ge0$ and $V_n\le V_{n+1}$, we have $$ V_nu_n-V_{n+1}u_{n+1} \le V_{n+1}(u_n-u_{n+1}) = V_{n+1}z. $$ 

Therefore, by \eqref{eq:Vn} and the Hardy inequality \eqref{eq:hardy}, $$ \begin{aligned} \mu\int_\Omega \bigl(V_nu_n-V_{n+1}u_{n+1}\bigr)z^+\dd x \le \mu\int_\Omega V_{n+1}(z^+)^2\dd x \le \frac{\mu}{\LN}\|\nabla z^+\|_2^2 \le \frac{\mu}{\LN}\rho(z^+)^2. \end{aligned} $$

Finally, since $g_n\le g_{n+1}$ and $$ \frac{1}{(u_n+1/n)^\gamma} \le \frac{1}{(u_n+1/(n+1))^\gamma}, $$ we have $$ \begin{aligned} &\frac{g_n}{(u_n+1/n)^\gamma} -\frac{g_{n+1}}{(u_{n+1}+1/(n+1))^\gamma} \\ &\qquad\le g_{n+1} \left( \frac{1}{(u_n+1/(n+1))^\gamma} -\frac{1}{(u_{n+1}+1/(n+1))^\gamma} \right). \end{aligned} $$ 

The right-hand side is nonpositive on $ \{z>0\}=\{u_n>u_{n+1}\}, $ because the function $ t\mapsto\frac{1}{(t+1/(n+1))^\gamma} $ is decreasing. Hence the right-hand side of \eqref{delta1} is nonpositive. Combining the preceding estimates, we obtain $$ \left(1-\frac{\mu}{\LN}\right)\rho(z^+)^2\le0. $$ Since $\mu<\LN$, it follows that $\rho(z^+)=0$, and hence $z^+=0$ a.e. in $\Omega$. Therefore, $ u_n\le u_{n+1}$ a.e. in $\Omega. $

We now prove the second assertion. By Lemma~\ref{lem:Pn}, $u_1\in H^1_0(\Omega)$ and $u_1\ge0$ a.e. in $\RN$. Moreover, by \eqref{weak11}, $$ \cE(u_1,\varphi) = \int_\Omega \left( \mu V_1u_1+\frac{g_1}{(u_1+1)^\gamma} \right)\varphi\dd x \ge0 $$
for every nonnegative $\varphi\in C^\infty_c(\Omega)$. Thus $u_1$ is a nonnegative weak supersolution of $ \cL u_1=0 \quad\text{in }\Omega. $ 

Furthermore, $u_1\not\equiv0$ in $\Omega$. Indeed, if $u_1\equiv0$, then \eqref{weak11} would give $$ \int_\Omega g_1\psi\dd x=0 $$ for every $\psi\in H^1_0(\Omega)$, and hence for every $\psi\in C^\infty_c(\Omega)$. Therefore, $g_1=0$ a.e. in $\Omega$. Since $g_1=\min\{g,1\}$ and $g\ge0$, this would imply $g=0$ a.e. in $\Omega$, contrary to the assumption $g\not\equiv0$. Lemma~\ref{lem:strongmin} now implies that, for every $\omega\Subset\Omega$, there exists $c_\omega>0$ such that $ u_1\ge c_\omega$ a.e. in $\omega. $ By the monotonicity of $(u_n)$, $ u_n\ge u_1\ge c_\omega$ a.e. in  $\omega $ for every $n\in\N$. This proves \eqref{eq:unifpos}. \end{proof}

\begin{remark}\label{rem:limit} By Lemma~\ref{lem:monotone}, the pointwise limit $$ u(x):=\lim_{n\to\infty}u_n(x) =\sup_{n\in\N}u_n(x)\in[0,\infty] $$ exists for a.e. $x\in\RN$. Moreover, $u=0$ a.e. in $\RN\setminus\Omega$ and, for every $\omega\Subset\Omega$, there exists $c_\omega>0$ such that $ u\ge c_\omega$ a.e. in $\omega. $ If $(u_n)$ is bounded in $L^1(\Omega)$, then the monotone convergence theorem yields $u\in L^1(\Omega)$, and hence $u<\infty$ a.e. in $\Omega$. Furthermore, $$ \|u-u_n\|_{L^1(\Omega)} = \int_\Omega (u-u_n)\dd x \longrightarrow0, $$ so that $u_n\to u$ in $L^1(\Omega)$. This function $u$ will be shown to be a solution in the proof of Theorems~\ref{thm:main1}--\ref{thm:main3}. \end{remark}

\section{A priori estimates}\label{sec:estimates}
Throughout this section, $(u_n)$ denotes the sequence of weak solutions of \eqref{problem11} given by Lemma~\ref{lem:Pn}. The estimates follow the strategy developed in \cite[Sections 3--5]{boccardo2010semilinear}. The main additional point in the present setting is the absorption of the Hardy term, which leads to the restrictions imposed on $\mu$. We recall that $$ u_n\in H^1_0(\Omega)\cap L^\infty(\Omega), \qquad u_n>0 \quad\text{a.e. in }\Omega, \qquad u_n\le u_{n+1}. $$ Moreover, by \eqref{eq:unifpos}, the sequence $(u_n)$ is uniformly bounded away from zero on every compact subset of $\Omega$. Throughout this section, $C$ denotes a positive constant independent of $n$, whose value may change from line to line.

\subsection{\texorpdfstring{The infinite-energy case $0\le\gamma<1$, $m<m^\#$}{The infinite-energy case}}

\begin{lemma}\label{prioriestimates1}
Let $0\le\gamma<1$, $(m,\gamma)\in\cS$ with $1\le m<m^\#$, let $g\in L^m(\Omega)$ be nonnegative and not identically zero, and let $0<\mu<\mu_{m,\gamma}$. Set $\delta:=\delta_m\in(0,1)$ and $v_n:=u_n^{(1+\delta)/2}$. Then:
\begin{enumerate}[label=(\roman*)]
\item $v_n\in H^1_0(\Omega)$, $\nabla v_n=\frac{1+\delta}{2}\,u_n^{\frac{\delta-1}{2}}\nabla u_n$ a.e. in $\Omega$, and $\rho(v_n)\le C$;
\item $\|u_n\|_{\sigma_m}\le C$;
\item $\|\nabla u_n\|_{q_m}\le C$; in particular, $(u_n)$ is bounded in $W^{1,q_m}_0(\Omega)$.
\end{enumerate}
\end{lemma}

\begin{proof} Throughout the proof, we use the constant \begin{equation}\label{eq:Kdelta} K_{\delta,\mu} := \frac{4\delta}{(1+\delta)^2} - \frac{\mu}{\LN}>0. \end{equation} 

Its positivity follows from the assumption $\mu<\mu_{m,\gamma}$ and the definition \eqref{eq:threshold}. 

 For $\varepsilon\in(0,1)$, define $ \Phi_\varepsilon(t) := (t+\varepsilon)^\delta-\varepsilon^\delta, \,\, t\ge0, $ and $ W_\varepsilon := (u_n+\varepsilon)^{\frac{1+\delta}{2}} - \varepsilon^{\frac{1+\delta}{2}}. $ Since $0<\delta<1$, the functions $ t\mapsto\Phi_\varepsilon(t) $ and $ t\longmapsto (t+\varepsilon)^{\frac{1+\delta}{2}} -\varepsilon^{\frac{1+\delta}{2}} $ are Lipschitz on $[0,\infty)$ and vanish at $t=0$. Therefore, $ \Phi_\varepsilon(u_n),W_\varepsilon\in H^1_0(\Omega). $
 
 Moreover, $$ \begin{aligned} \nabla u_n\cdot\nabla\Phi_\varepsilon(u_n) = \delta(u_n+\varepsilon)^{\delta-1}|\nabla u_n|^2 = \frac{4\delta}{(1+\delta)^2} |\nabla W_\varepsilon|^2 \end{aligned} $$ a.e. in $\Omega$. Applying \eqref{eq:alg} with $\theta=\delta$, $a=u_n(x)+\varepsilon$, and $b=u_n(y)+\varepsilon$, we also obtain $$ \begin{aligned} (u_n(x)-u_n(y)) \bigl( \Phi_\varepsilon(u_n(x)) - \Phi_\varepsilon(u_n(y)) \bigr) \ge \frac{4\delta}{(1+\delta)^2} \bigl( W_\varepsilon(x)-W_\varepsilon(y) \bigr)^2 \end{aligned} $$ for every $x,y\in\RN$. Consequently, $$ \cE\bigl(u_n,\Phi_\varepsilon(u_n)\bigr) \ge \frac{4\delta}{(1+\delta)^2} \rho(W_\varepsilon)^2. $$

 We may use $\Phi_\varepsilon(u_n)$ as a test function in \eqref{weak11}. Since $0<\delta<1$, the subadditivity of the map $t\mapsto t^\delta$ gives $ 0\le\Phi_\varepsilon(t)\le t^\delta$ for every $t\ge0. $ 
 
 It follows that \begin{equation}\label{eq:regenergy} \begin{aligned} \frac{4\delta}{(1+\delta)^2} \rho(W_\varepsilon)^2 &\le \mu\int_\Omega V_nu_n\Phi_\varepsilon(u_n)\dd x + \int_\Omega \frac{g_n\Phi_\varepsilon(u_n)} {(u_n+1/n)^\gamma}\dd x \\ &\le \mu\int_\Omega V_nu_n^{1+\delta}\dd x + \int_\Omega g_nu_n^{\delta-\gamma}\dd x. \end{aligned} \end{equation} In the last inequality, we used $ (u_n+1/n)^\gamma\ge u_n^\gamma $ and the fact that $u_n>0$ a.e. in $\Omega$. Notice also that $\delta-\gamma\ge0$ by Lemma~\ref{lem:exponents}\textup{(ii)}.  For each fixed $n$, the right-hand side of \eqref{eq:regenergy} is finite and independent of $\varepsilon$, because $u_n\in L^\infty(\Omega)$ and $V_n\le n$.

 We now let $\varepsilon\to0$. By \eqref{eq:regenergy} and \eqref{eq:normequiv}, the family $(W_\varepsilon)$ is bounded in $H^1_0(\Omega)$.  Moreover, $ W_\varepsilon\to v_n$ a.e. in $\Omega $ as $\varepsilon\to0$, and $ 0\le W_\varepsilon \le (u_n+1)^{\frac{1+\delta}{2}}. $
 Since $u_n\in L^\infty(\Omega)$ and $\Omega$ is bounded, the function on the right-hand side belongs to $L^2(\Omega)$. 
 Hence, by the dominated convergence theorem, $ W_\varepsilon\to v_n \quad\text{in }L^2(\Omega). $ Therefore, along a sequence $\varepsilon_j\to0$, $ W_{\varepsilon_j}\rightharpoonup v_n$ weakly in $H^1_0(\Omega). $ In particular, $v_n\in H^1_0(\Omega)$. 
 
 Since $\rho$ is a Hilbert norm on $H^1_0(\Omega)$, its weak lower semicontinuity gives $$ \rho(v_n)^2 \le \liminf_{j\to\infty} \rho(W_{\varepsilon_j})^2. $$ Taking the lower limit in \eqref{eq:regenergy}, using $u_n^{1+\delta}=v_n^2$, and then applying \eqref{eq:Vn} and the Hardy inequality \eqref{eq:hardy}, we obtain $$ \begin{aligned} \frac{4\delta}{(1+\delta)^2}\rho(v_n)^2 &\le \mu\int_\Omega V_nv_n^2\dd x + \int_\Omega g_nu_n^{\delta-\gamma}\dd x \le \mu\int_\Omega\frac{v_n^2}{|x|^2}\dd x + \int_\Omega g_nu_n^{\delta-\gamma}\dd x \\ &\le \frac{\mu}{\LN}\rho(v_n)^2 + \int_\Omega g_nu_n^{\delta-\gamma}\dd x. \end{aligned} $$ Thus, \begin{equation}\label{kdeltar} K_{\delta,\mu}\rho(v_n)^2 \le \int_\Omega g_nu_n^{\delta-\gamma}\dd x. \end{equation}

 We next establish the formula for $\nabla v_n$. For $t>0$, define $ \Psi_t := \bigl(\max\{u_n,t\}\bigr)^{\frac{1+\delta}{2}} -t^{\frac{1+\delta}{2}}. $  The function $ r\longmapsto \bigl(\max\{r,t\}\bigr)^{\frac{1+\delta}{2}} - t^{\frac{1+\delta}{2}} $ is Lipschitz on $[0,\infty)$ and vanishes at $r=0$. Hence $\Psi_t\in H^1_0(\Omega)$, and the  chain rule gives $$ \nabla\Psi_t = \frac{1+\delta}{2} u_n^{\frac{\delta-1}{2}}\nabla u_n \chi_{\{u_n>t\}} $$ a.e. in $\Omega$. On the set $\{u_n>t\}$, we have $ v_n-\Psi_t=t^{\frac{1+\delta}{2}}. $ Therefore, $$ \nabla v_n = \nabla\Psi_t = \frac{1+\delta}{2} u_n^{\frac{\delta-1}{2}}\nabla u_n $$ a.e. on $\{u_n>t\}$. Taking $t=1/k$ and letting $k\to\infty$, and using $u_n>0$ a.e. in $\Omega$, we conclude that $$ \nabla v_n = \frac{1+\delta}{2} u_n^{\frac{\delta-1}{2}}\nabla u_n $$ a.e. in $\Omega$. 
 
We now derive the estimates that are uniform in $n$. If $m=1$, then $\delta=\gamma$, and \eqref{kdeltar} gives $ K_{\delta,\mu}\rho(v_n)^2\le\|g\|_1. $ Thus, $ \rho(v_n)\le C. $

Suppose now that $m>1$. By H\"older's inequality, Lemma~\ref{lem:exponents}\textup{(iii)}, and the Sobolev inequality, we obtain $$ \begin{aligned} \int_\Omega g_nu_n^{\delta-\gamma}\dd x \le \|g\|_m \left( \int_\Omega u_n^{(\delta-\gamma)m'}\dd x \right)^{1/m'} = \|g\|_m \left( \int_\Omega v_n^{2^*}\dd x \right)^{1/m'} \le \|g\|_m \left( S_N\rho(v_n)^2 \right)^{\frac{2^*}{2m'}}. \end{aligned} $$ 

Combining this estimate with \eqref{kdeltar}, we find $$ \rho(v_n)^{2-\frac{2^*}{m'}} \le \frac{ S_N^{\frac{2^*}{2m'}} \|g\|_m }{ K_{\delta,\mu} }. $$ Moreover, $ 2-\frac{2^*}{m'} = \frac{2(N-2m)}{m(N-2)}>0, $ because $m<N/2$. Hence, $ \rho(v_n)\le C. $ This proves the uniform bound in \textup{(i)}. By Lemma~\ref{lem:exponents}\textup{(iii)}, $ u_n^{\sigma_m}=v_n^{2^*}. $ For $m=1$, this identity follows directly from $ \sigma_1=\frac{N(1+\gamma)}{N-2} $ and $\delta=\gamma$. Therefore, the Sobolev inequality gives $$ \begin{aligned} \int_\Omega u_n^{\sigma_m}\dd x = \int_\Omega v_n^{2^*}\dd x \le \left( S_N\rho(v_n)^2 \right)^{2^*/2} \le C. \end{aligned} $$ Consequently, $ \|u_n\|_{\sigma_m}\le C, $ which proves \textup{(ii)}.

It remains to prove \textup{(iii)}. Set $q:=q_m\in(1,2)$. By the gradient identity established above, $$ u_n^{\delta-1}|\nabla u_n|^2 = \frac{4}{(1+\delta)^2}|\nabla v_n|^2 $$ a.e. in $\Omega.$ Hence, $ |\nabla u_n|^q = \left( u_n^{\delta-1}|\nabla u_n|^2 \right)^{q/2} u_n^{\frac{(1-\delta)q}{2}} $ a.e. in $\Omega$. Applying H\"older's inequality with conjugate exponents $2/q$ and $2/(2-q)$, we obtain \begin{equation}\label{nablaun} \begin{aligned} \int_\Omega|\nabla u_n|^q\dd x &\le \left( \frac{4}{(1+\delta)^2} \int_\Omega|\nabla v_n|^2\dd x \right)^{q/2}  \left( \int_\Omega u_n^{\frac{(1-\delta)q}{2-q}}\dd x \right)^{\frac{2-q}{2}}. \end{aligned} \end{equation} By Lemma~\ref{lem:exponents}\textup{(iv)}, $ \frac{(1-\delta)q}{2-q}=\sigma_m. $ Thus, \eqref{nablaun} becomes $$ \int_\Omega|\nabla u_n|^q\dd x \le \left( \frac{4\|\nabla v_n\|_2^2}{(1+\delta)^2} \right)^{q/2} \left( \int_\Omega u_n^{\sigma_m}\dd x \right)^{\frac{2-q}{2}}. $$ Both factors on the right-hand side are uniformly bounded by \textup{(i)} and \textup{(ii)}. Therefore, $ \|\nabla u_n\|_{q_m}\le C. $ Finally, since $u_n\in H^1_0(\Omega)\subset W^{1,q_m}_0(\Omega)$, the Poincar\'e inequality yields $ \|u_n\|_{W^{1,q_m}_0(\Omega)}\le C. $ This proves \textup{(iii)}. \end{proof}

\subsection{The finite-energy case}

\begin{lemma}\label{prioriestimates2}
Let $(m,\gamma)\in\cS$ with either $\gamma=1$, or $0\le\gamma<1$ and $m\ge m^\#$, let $g\in L^m(\Omega)$ be nonnegative and not identically zero, and let $0<\mu<\LN$. Then $(u_n)$ is bounded in $H^1_0(\Omega)$.
\end{lemma}

\begin{proof}
Testing \eqref{weak11} with $\psi=u_n$ and using \eqref{coercive} we get
\begin{equation}\label{eq:energy-un}
\Big(1-\frac{\mu}{\LN}\Big)\rho(u_n)^2\le\int_\Omega\frac{g_n\,u_n}{(u_n+1/n)^\gamma}\dd x\le\int_\Omega g\,u_n^{1-\gamma}\dd x .
\end{equation}
If $\gamma=1$ the right-hand side is bounded by $\|g\|_1$. If $0\le\gamma<1$ and $m\ge m^\#$, then $(1-\gamma)m'\le2^*$ (with $m'=1$ if $m=\infty$); by H\"older's inequality, the boundedness of $\Omega$ and \eqref{eq:sobolev},
\[
\int_\Omega g\,u_n^{1-\gamma}\dd x\le\|g\|_m\|u_n\|_{(1-\gamma)m'}^{1-\gamma}\le C\|g\|_m\|u_n\|_{2^*}^{1-\gamma}\le C\|g\|_m\,\rho(u_n)^{1-\gamma}.
\]
Since $1-\gamma<2$, \eqref{eq:energy-un} implies $\rho(u_n)\le C$.
\end{proof}

\subsection{\texorpdfstring{The strongly singular case $\gamma>1$}{The strongly singular case}}

\begin{lemma}\label{prioriestimates3}
Let $\gamma>1$, let $g\in L^1(\Omega)$ be nonnegative and not identically zero, and let $0<\mu<\frac{4\gamma}{(1+\gamma)^2}\LN$. Then the following assertions hold:
\begin{enumerate}[label=(\roman*)]
\item $(u_n^{(\gamma+1)/2})$ is bounded in $H^1_0(\Omega).$ In particular $(u_n)$ is bounded in $L^{r}(\Omega)$ with  $r:=\frac{(\gamma+1)2^*}{2}>2.$
\item For every $\varphi\in C^\infty_c(\Omega)$ with $0\le\varphi\le1$, the sequence $(u_n\varphi)$ is bounded in $H^1_0(\Omega).$ Consequently, $(u_n)$ is bounded in $H^1(\omega)$ for every $\omega\Subset\Omega$.
\end{enumerate}
\end{lemma}

\begin{proof}
(i) Since $u_n\in L^\infty(\Omega)$ and $\gamma>1$, the functions $t\mapsto t^\gamma$ and $t\mapsto t^{(\gamma+1)/2}$ are Lipschitz on $[0,\|u_n\|_\infty]$ and vanish at $0$; hence $u_n^\gamma,\ w_n:=u_n^{(\gamma+1)/2}\in H^1_0(\Omega)$ and by chain rule, we get
\[
\nabla u_n\cdot\nabla u_n^\gamma=\gamma u_n^{\gamma-1}|\nabla u_n|^2=\frac{4\gamma}{(\gamma+1)^2}|\nabla w_n|^2\qquad\text{a.e. in }\Omega .
\]

Moreover, applying \eqref{eq:alg} with $\theta=\gamma$, we obtain $$ \begin{aligned} (u_n(x)-u_n(y)) \bigl(u_n(x)^\gamma-u_n(y)^\gamma\bigr) \ge \frac{4\gamma}{(\gamma+1)^2} \bigl(w_n(x)-w_n(y)\bigr)^2 \end{aligned} $$ for every $x,y\in\RN$. It follows that $$ \cE(u_n,u_n^\gamma) \ge \frac{4\gamma}{(\gamma+1)^2}\rho(w_n)^2. $$

 Testing \eqref{weak11} with $u_n^\gamma$ and using \eqref{eq:Vn} and the Hardy inequality \eqref{eq:hardy}, we obtain
\[
\frac{4\gamma}{(\gamma+1)^2}\rho(w_n)^2\le\mu\int_\Omega\frac{w_n^2}{|x|^2}\dd x+\int_\Omega g_n\frac{u_n^\gamma}{(u_n+1/n)^\gamma}\dd x\le\frac{\mu}{\LN}\rho(w_n)^2+\|g\|_1,
\]
so that
\[
\Big(\frac{4\gamma}{(\gamma+1)^2}-\frac{\mu}{\LN}\Big)\rho(w_n)^2\le\|g\|_1 ,
\]
and the coefficient on the left-hand side is positive by assumption. Hence, $ \rho(w_n)\le C, $ and $(w_n)$ is bounded in $H^1_0(\Omega)$.  Set $$ r:=\frac{(\gamma+1)2^*}{2}. $$ Since $ u_n^r=w_n^{2^*}, $ the Sobolev inequality gives $$ \begin{aligned} \int_\Omega u_n^r\dd x = \int_\Omega w_n^{2^*}\dd x \le \left(S_N\rho(w_n)^2\right)^{2^*/2} \le C. \end{aligned} $$ Thus, $(u_n)$ is bounded in $L^r(\Omega)$. Since $\gamma>1$, we have $r>2^*>2$, and therefore $(u_n)$ is also bounded in $L^2(\Omega)$. This proves \textup{(i)}.

(ii) Let $\varphi\in C^\infty_c(\Omega)$, $0\le\varphi\le1$, and $K:=\supp\varphi$. Since $u_n\varphi^2\in H^1_0(\Omega)$, we can test \eqref{weak11} with it. For the local term we have the identity
\[
\int_\Omega\nabla u_n\cdot\nabla(u_n\varphi^2)\dd x=\int_\Omega|\nabla(u_n\varphi)|^2\dd x-\int_\Omega u_n^2|\nabla\varphi|^2\dd x .
\]
For the nonlocal term we use the identity
\[
(a-b)(aA^2-bB^2)=(aA-bB)^2-ab(A-B)^2,\qquad a,b,A,B\in\R,
\]
with $a=u_n(x)$, $b=u_n(y)$, $A=\varphi(x)$, $B=\varphi(y)$, together with $ab\le\frac12(a^2+b^2)$ and the symmetry of the kernel, to get
\begin{multline*}
\iint_{\R^{2N}}\frac{(u_n(x)-u_n(y))\big((u_n\varphi^2)(x)-(u_n\varphi^2)(y)\big)}{|x-y|^{N+2s}}\dd x\dd y\\
\ge[u_n\varphi]_s^2-\int_{\RN}u_n(x)^2\int_{\RN}\frac{(\varphi(x)-\varphi(y))^2}{|x-y|^{N+2s}}\dd y\dd x .
\end{multline*}
Since $|\varphi(x)-\varphi(y)|\le\min\{\|\nabla\varphi\|_\infty|x-y|,1\}$ and $s \in (0, 1)$,
\[
\sup_{x\in\RN}\int_{\RN}\frac{(\varphi(x)-\varphi(y))^2}{|x-y|^{N+2s}}\dd y\le\|\nabla\varphi\|_\infty^2\int_{\{|z|<1\}}\frac{\dd z}{|z|^{N+2s-2}}+\int_{\{|z|\ge1\}}\frac{\dd z}{|z|^{N+2s}}=:M_\varphi<\infty .
\]
For the Hardy term, \eqref{eq:Vn} and the Hardy inequality applied to $u_n\varphi\in H^1_0(\Omega)$ give
\[
\mu\int_\Omega V_nu_n^2\varphi^2\dd x\le\frac{\mu}{\LN}\int_\Omega|\nabla(u_n\varphi)|^2\dd x .
\]
Finally, by \eqref{eq:unifpos}, there exists $c_K>0$, independent of $n$, such that $ u_n\ge c_K$ a.e. in $K.$  Since $\gamma>1$, we have
\[
\int_\Omega\frac{g_n\,u_n\varphi^2}{(u_n+1/n)^\gamma}\dd x\le\int_K g\,u_n^{1-\gamma}\dd x\le c_K^{1-\gamma}\|g\|_1 .
\]

Combining the preceding estimates in the weak formulation, we find $$ \begin{aligned} &\left(1-\frac{\mu}{\LN}\right) \int_\Omega|\nabla(u_n\varphi)|^2\dd x + \frac{\Cs}{2}[u_n\varphi]_s^2 \le \left( \|\nabla\varphi\|_\infty^2 + \frac{\Cs}{2}M_\varphi \right) \|u_n\|_2^2 + c_K^{1-\gamma}\|g\|_1. \end{aligned} $$

Since the nonlocal term on the left-hand side is nonnegative, we may discard it and obtain $$ \left(1-\frac{\mu}{\LN}\right) \int_\Omega|\nabla(u_n\varphi)|^2\dd x \le \left( \|\nabla\varphi\|_\infty^2 + \frac{\Cs}{2}M_\varphi \right) \|u_n\|_2^2 + c_K^{1-\gamma}\|g\|_1. $$ The right-hand side is uniformly bounded by \textup{(i)}. Since $\mu<\frac{4\gamma}{(1+\gamma)^2}\LN<\LN$, the coefficient on the left-hand side is positive. Thus, $ \|\nabla(u_n\varphi)\|_2\le C_\varphi. $ The Poincar\'e inequality now implies that $(u_n\varphi)$ is bounded in $H^1_0(\Omega)$. Given $\omega\Subset\Omega$, choose $\varphi\in C^\infty_c(\Omega)$ such that $0\le\varphi\le1$ and $\varphi\equiv1$ in a neighborhood of $\overline\omega$. Then $ u_n\varphi=u_n$ in $\omega, $ and hence $$ \|u_n\|_{H^1(\omega)} \le \|u_n\varphi\|_{H^1(\Omega)} \le C_\omega. $$ This proves \textup{(ii)}.
\end{proof}

\begin{remark}\label{rem:local-mu} In Lemma~\ref{prioriestimates3}, the restriction $ \mu< \frac{4\gamma}{(1+\gamma)^2}\LN $ is used only to obtain the global estimate in \textup{(i)}. Once a global bound for $(u_n)$ in $L^2(\Omega)$ is available, the local estimate in \textup{(ii)} only requires $ \mu<\LN. $ \end{remark}

\section{Proofs of the main results}\label{sec:proofs}

We first isolate the passage-to-the-limit argument, which is common to the proofs of the three theorems. \begin{proposition}\label{prop:limit} Let $g\in L^1(\Omega)$ be nonnegative and not identically zero, let $\gamma\ge0$ and $0<\mu<\LN$, and let $(u_n)$ be the sequence given by Lemma~\ref{lem:Pn}. Assume that there exists $ q\in\left(\frac{N}{N-1},N\right) $ such that $(u_n)$ is bounded in $L^1(\Omega)$ and, for every open set $U\Subset\Omega$, bounded in $W^{1,q}(U)$. Then the function $ u:=\lim_{n\to\infty}u_n $ introduced in Remark~\ref{rem:limit} is a weak solution of \eqref{problem1} in the sense of Definition~\ref{def:weak}. Moreover, for every open set $U\Subset\Omega$, $$ u_n\rightharpoonup u \quad\text{weakly in }W^{1,q}(U). $$ \end{proposition} 

\begin{proof} By Remark~\ref{rem:limit}, we have $ u\in L^1(\Omega)$ with $u=0$ a.e. in $\RN\setminus\Omega, $ and $ u_n\to u $ in $L^1(\Omega). $ Moreover, for every $\omega\Subset\Omega$, there exists $c_\omega>0$ such that $$ u\ge u_n\ge u_1\ge c_\omega \quad\text{a.e. in }\omega. $$ 

Thus, $u$ satisfies condition \textup{(b)} of Definition~\ref{def:weak}. Fix an open set $U\Subset\Omega$. Since $(u_n)$ is bounded in $W^{1,q}(U)$ and $1<q<\infty$, every subsequence of $(u_n)$ admits a further subsequence converging weakly in $W^{1,q}(U)$ to some function $v\in W^{1,q}(U)$. On the other hand, $$ u_n\to u \quad\text{in }L^1(U), $$ so the same further subsequence converges to $u$ in $L^1(U)$. Hence $v=u$ a.e. in $U$. Since every subsequence admits a further subsequence converging weakly to $u$, the whole sequence satisfies $$ u_n\rightharpoonup u \quad\text{weakly in }W^{1,q}(U). $$ 

In particular, $ u\in W^{1,q}(U). $ Since $U\Subset\Omega$ is arbitrary, it follows that $ u\in W^{1,q}_{\rm loc}(\Omega), $ and therefore $u$ satisfies condition \textup{(a)} of Definition~\ref{def:weak}. It remains to pass to the limit in the weak formulation. Fix $\psi\in C^\infty_c(\Omega)$ and set $ K:=\supp\psi. $ By Lemma~\ref{lem:welldef}\textup{(iii)}, $ \cE(u_n,\psi)\to\cE(u,\psi). $ 

We next consider the Hardy term. Write $$ \psi=\psi^+-\psi^-, \qquad \psi^+:=\max\{\psi,0\}, \qquad \psi^-:=\max\{-\psi,0\}. $$ Since both $(V_n)$ and $(u_n)$ are nondecreasing and nonnegative, the sequences $ V_nu_n\psi^+$ and $V_nu_n\psi^- $ are nondecreasing. Moreover, $$ V_nu_n\psi^\pm \longrightarrow \frac{u\psi^\pm}{|x|^2} \quad\text{a.e. in }\Omega. $$ By Lemma~\ref{lem:welldef}\textup{(i)}, $ \frac{u|\psi|}{|x|^2}\in L^1(\Omega). $ The monotone convergence theorem therefore gives $$ \int_\Omega V_nu_n\psi^\pm\dd x \longrightarrow \int_\Omega\frac{u\psi^\pm}{|x|^2}\dd x. $$ Subtracting the two limits, we obtain $$ \int_\Omega V_nu_n\psi\dd x \longrightarrow \int_\Omega\frac{u\psi}{|x|^2}\dd x. $$ 

Finally, we consider the singular term. By \eqref{eq:unifpos}, there exists $c_K>0$, independent of $n$, such that $ u_n\ge c_K $ a.e. in $K. $ Consequently, $$ \left| \frac{g_n\psi}{(u_n+1/n)^\gamma} \right| \le c_K^{-\gamma}\|\psi\|_\infty g\chi_K \quad\text{a.e. in }\Omega. $$

The right-hand side belongs to $L^1(\Omega)$. Furthermore, $ g_n\to g$ and $u_n+\frac1n\to u $ a.e. in $\Omega$. Since $u\ge c_K>0$ a.e. on $K$, it follows that $$ \frac{g_n\psi}{(u_n+1/n)^\gamma} \longrightarrow \frac{g\psi}{u^\gamma} \quad\text{a.e. in }\Omega. $$ Thus, the dominated convergence theorem yields $$ \int_\Omega \frac{g_n\psi}{(u_n+1/n)^\gamma}\dd x \longrightarrow \int_\Omega\frac{g\psi}{u^\gamma}\dd x. $$ 

Since $\psi\in C^\infty_c(\Omega)\subset H^1_0(\Omega)$, it is an admissible test function in \eqref{weak11}. Passing to the limit in that identity, we obtain $$ \cE(u,\psi) - \mu\int_\Omega\frac{u\psi}{|x|^2}\dd x = \int_\Omega\frac{g\psi}{u^\gamma}\dd x. $$ 

Therefore, $u$ satisfies \eqref{eq:weakform} for every $\psi\in C^\infty_c(\Omega)$ and is a weak solution of \eqref{problem1}. \end{proof}

\begin{proof}[Proof of Theorem~\ref{thm:main1}] By Lemma~\ref{prioriestimates1}, the sequence $(u_n)$ is bounded in $W^{1,q_m}_0(\Omega)$. Moreover, by Lemma~\ref{lem:exponents}\textup{(iv)}, $ \frac{N}{N-1}<q_m<2<N. $ In particular, $(u_n)$ is bounded in $L^1(\Omega)$ and in $W^{1,q_m}(U)$ for every open set $U\Subset\Omega$. Therefore, Proposition~\ref{prop:limit} shows that the function $ u=\lim_{n\to\infty}u_n $ is a weak solution of \eqref{problem1}. Since $W^{1,q_m}_0(\Omega)$ is reflexive, every subsequence of $(u_n)$ admits a further subsequence converging weakly in $W^{1,q_m}_0(\Omega)$. On the other hand, by Remark~\ref{rem:limit}, $ u_n\to u \quad\text{in }L^1(\Omega). $ Thus, every such weak limit coincides with $u$. Consequently, $ u\in W^{1,q_m}_0(\Omega) $ and the whole sequence satisfies $$ u_n\rightharpoonup u \quad\text{weakly in }W^{1,q_m}_0(\Omega). $$

Next, set $ v_n:=u_n^{\frac{1+\delta_m}{2}}. $ By Lemma~\ref{prioriestimates1}\textup{(i)}, $(v_n)$ is bounded in $H^1_0(\Omega)$. Hence, every subsequence of $(v_n)$ admits a further subsequence, still denoted by $(v_n)$, and a function $v\in H^1_0(\Omega)$ such that $ v_n\rightharpoonup v $ weakly in $H^1_0(\Omega). $ By the Rellich--Kondrachov theorem, after passing to a further subsequence, $ v_n\to v $ strongly in $L^2(\Omega) $ and a.e. in $\Omega$. Since $u_n\to u$ a.e. in $\Omega$, we also have $ v_n = u_n^{\frac{1+\delta_m}{2}} \longrightarrow u^{\frac{1+\delta_m}{2}}$ a.e. in  $\Omega. $ Therefore, $ v=u^{\frac{1+\delta_m}{2}}$ a.e. in  $\Omega, $ and hence $ u^{\frac{1+\delta_m}{2}}\in H^1_0(\Omega). $ Finally, Lemma~\ref{prioriestimates1}\textup{(ii)} gives $ \|u_n\|_{\sigma_m}\le C. $ Since $u_n\to u$ a.e. in $\Omega$, Fatou's lemma yields $$ \int_\Omega u^{\sigma_m}\dd x \le \liminf_{n\to\infty} \int_\Omega u_n^{\sigma_m}\dd x \le C. $$ Thus, $ u\in L^{\sigma_m}(\Omega). $ This completes the proof. \end{proof}

\begin{proof}[Proof of Theorem~\ref{thm:main2}] By Lemma~\ref{prioriestimates2}, $(u_n)$ is bounded in $H^1_0(\Omega)$. Since $\Omega$ is bounded, the Poincar\'e inequality implies that $(u_n)$ is bounded in $L^2(\Omega)$, and hence in $L^1(\Omega)$. Moreover, $ \frac{N}{N-1}<2<N $ because $N\ge3$. Therefore, Proposition~\ref{prop:limit}, applied with $q=2$, shows that $u$ is a weak solution of \eqref{problem1} and that $ u_n\rightharpoonup u$ weakly in $H^1(U) $ for every open set $U\Subset\Omega$. Since $(u_n)$ is bounded in the reflexive space $H^1_0(\Omega)$, every subsequence admits a further subsequence converging weakly in $H^1_0(\Omega)$ to some function $v\in H^1_0(\Omega)$. On the other hand, by Remark~\ref{rem:limit}, $ u_n\to u$ in  $L^1(\Omega). $ Hence $v=u$ a.e. in $\Omega$. It follows that $u\in H^1_0(\Omega) $ and that the whole sequence satisfies $ u_n\rightharpoonup u $ weakly in $H^1_0(\Omega). $ This completes the proof. \end{proof}

\begin{proof}[Proof of Theorem~\ref{thm:main3}] By Lemma~\ref{prioriestimates3}\textup{(i)}, $(u_n)$ is bounded in $L^r(\Omega)$ for some $r>2$. Since $\Omega$ is bounded, it follows that $(u_n)$ is bounded in $L^2(\Omega)$, and hence in $L^1(\Omega)$. Moreover, by Lemma~\ref{prioriestimates3}\textup{(ii)}, $(u_n)$ is bounded in $H^1(U)$ for every open set $U\Subset\Omega$. Therefore, Proposition~\ref{prop:limit}, applied with $q=2$, shows that $ u\in H^1_{\rm loc}(\Omega) $ is a weak solution of \eqref{problem1}. It also gives $ u_n\rightharpoonup u$ weakly in $H^1(U) $ for every open set $U\Subset\Omega$. Finally, set $ w_n:=u_n^{\frac{\gamma+1}{2}}. $ By Lemma~\ref{prioriestimates3}\textup{(i)}, $(w_n)$ is bounded in $H^1_0(\Omega)$. Thus, every subsequence of $(w_n)$ admits a further subsequence and a function $w\in H^1_0(\Omega)$ such that $ w_n\rightharpoonup w$ weakly in $H^1_0(\Omega). $ By the Rellich--Kondrachov theorem, after passing to a further subsequence, $ w_n\to w$ strongly in $L^2(\Omega) $ and a.e. in $\Omega$. On the other hand, since $u_n\to u$ a.e. in $\Omega$, we have $ w_n = u_n^{\frac{\gamma+1}{2}} \longrightarrow u^{\frac{\gamma+1}{2}}$ a.e. in $\Omega. $ Consequently, $ w=u^{\frac{\gamma+1}{2}}$ a.e. in $\Omega, $ and hence $u^{\frac{\gamma+1}{2}}\in H^1_0(\Omega). $ This completes the proof. \end{proof}

\begin{remark}\label{rem:final}
\begin{enumerate}[label=(\roman*), leftmargin=*]
\item The solutions constructed above are limits of the nondecreasing sequence $(u_n)$, and they satisfy $u\ge u_1$; in particular, they are bounded away from zero on compact subsets of $\Omega$, uniformly with respect to the data $g_n$.
\item Uniqueness in the energy space is proved in Theorem \ref{thm:uniqueness} below. Uniqueness of the infinite-energy solutions of Theorems \ref{thm:main1} and \ref{thm:main3}, as well as the optimality of the existence thresholds $\mu_{m,\gamma}$ (for $\delta_m<1$) and $\frac{4\gamma}{(1+\gamma)^2}\LN$, are left for future investigation.
\end{enumerate}
\end{remark}

\section{Summability, uniqueness and nonexistence}\label{sec:further}

\subsection{Proof of Theorem \ref{thm:summability}}
We first extend the estimate of Lemma \ref{prioriestimates1} to the case $\delta_m\ge1$.

\begin{lemma}\label{prioriestimates4}
Let $(m,\gamma)\in\cS$ with $1\le m<N/2$ and $\delta:=\delta_m\ge1$, let $g\in L^m(\Omega)$ be nonnegative and not identically zero, and let $0<\mu<\mu_{m,\gamma}$. Then $v_n:=u_n^{(1+\delta)/2}\in H^1_0(\Omega)$, $\rho(v_n)\le C$ and $\|u_n\|_{\sigma_m}\le C$, with $C$ independent of $n$.
\end{lemma}

\begin{proof} Since $u_n\in L^\infty(\Omega)$ and $\delta\ge1$, the functions $t\mapsto t^\delta$ and $t\mapsto t^{(1+\delta)/2}$ are Lipschitz on $[0,\|u_n\|_\infty]$ and vanish at $t=0$. Therefore, $ u_n^\delta,\, v_n:=u_n^{\frac{1+\delta}{2}} \in H^1_0(\Omega). $ Moreover, the chain rule gives $$ \begin{aligned} \nabla u_n\cdot\nabla u_n^\delta = \delta u_n^{\delta-1}|\nabla u_n|^2 = \frac{4\delta}{(1+\delta)^2}|\nabla v_n|^2 \end{aligned} $$ a.e. in $\Omega$. 

By \eqref{eq:alg}, applied with $\theta=\delta$, we also have $$ \begin{aligned} &(u_n(x)-u_n(y)) \bigl(u_n(x)^\delta-u_n(y)^\delta\bigr) \ge \frac{4\delta}{(1+\delta)^2} \bigl(v_n(x)-v_n(y)\bigr)^2 \end{aligned} $$ for every $x,y\in\RN$. Consequently, $$ \cE(u_n,u_n^\delta) \ge \frac{4\delta}{(1+\delta)^2}\rho(v_n)^2. $$ Testing \eqref{weak11} with $u_n^\delta$ and using \eqref{eq:Vn}, \eqref{eq:hardy}, and the fact that $u_n>0$ a.e. in $\Omega$, we obtain $$ \begin{aligned} \frac{4\delta}{(1+\delta)^2}\rho(v_n)^2 &\le \mu\int_\Omega V_nu_n^{1+\delta}\dd x + \int_\Omega \frac{g_nu_n^\delta}{(u_n+1/n)^\gamma}\dd x \\ &\le \mu\int_\Omega\frac{v_n^2}{|x|^2}\dd x + \int_\Omega g_nu_n^{\delta-\gamma}\dd x \le \frac{\mu}{\LN}\rho(v_n)^2 + \int_\Omega g_nu_n^{\delta-\gamma}\dd x. \end{aligned} $$ 

Here we have used $\delta-\gamma\ge0$, which follows from Lemma~\ref{lem:exponents}\textup{(ii)}. Hence, $$ K_{\delta,\mu}\rho(v_n)^2 \le \int_\Omega g_nu_n^{\delta-\gamma}\dd x, $$ where $ K_{\delta,\mu} := \frac{4\delta}{(1+\delta)^2} - \frac{\mu}{\LN}>0. $ 
The positivity of $K_{\delta,\mu}$ follows from $\mu<\mu_{m,\gamma}$ and \eqref{eq:threshold}. If $m=1$, then Lemma~\ref{lem:exponents}\textup{(i)} gives $\delta=\gamma$. Therefore, $$ K_{\delta,\mu}\rho(v_n)^2 \le \int_\Omega g_n\dd x \le \|g\|_1, $$ and hence $ \rho(v_n)\le C. $

Suppose now that $m>1$. By H\"older's inequality and Lemma~\ref{lem:exponents}\textup{(iii)}, $$ \begin{aligned} \int_\Omega g_nu_n^{\delta-\gamma}\dd x \le \|g\|_m \left( \int_\Omega u_n^{(\delta-\gamma)m'}\dd x \right)^{1/m'} = \|g\|_m \left( \int_\Omega v_n^{2^*}\dd x \right)^{1/m'}. \end{aligned} $$ 

The Sobolev inequality then yields $$ \int_\Omega g_nu_n^{\delta-\gamma}\dd x \le \|g\|_m \left( S_N\rho(v_n)^2 \right)^{\frac{2^*}{2m'}}. $$ Consequently, we have $$ \rho(v_n)^{2-\frac{2^*}{m'}} \le \frac{ S_N^{\frac{2^*}{2m'}}\|g\|_m }{ K_{\delta,\mu} }. $$

Since $m<N/2$, we have $ 2-\frac{2^*}{m'} = \frac{2(N-2m)}{m(N-2)}>0. $ It follows that $ \rho(v_n)\le C. $ Finally, Lemma~\ref{lem:exponents}\textup{(iii)} gives $ u_n^{\sigma_m}=v_n^{2^*}. $ Therefore, by the Sobolev inequality, $$ \begin{aligned} \int_\Omega u_n^{\sigma_m}\dd x = \int_\Omega v_n^{2^*}\dd x \le \left( S_N\rho(v_n)^2 \right)^{2^*/2} \le C. \end{aligned} $$ Thus, $ \|u_n\|_{\sigma_m}\le C, $ which completes the proof. \end{proof}

\begin{proof}[Proof of Theorem~\ref{thm:summability}] We first verify that one of Theorems~\ref{thm:main1}--\ref{thm:main3} provides a weak solution $ u=\lim_{n\to\infty}u_n $ of \eqref{problem1}. If $0\le\gamma<1$ and $m<m^\#$, then Theorem~\ref{thm:main1} applies, since its restriction on $\mu$ coincides with the present one. If $\gamma=1$, or if $0\le\gamma<1$ and $m\ge m^\#$, then Lemma~\ref{lem:exponents}\textup{(v)} gives $ \mu<\mu_{m,\gamma}\le\LN. $ Therefore, Theorem~\ref{thm:main2} applies. Finally, suppose that $\gamma>1$. By Lemma~\ref{lem:exponents}\textup{(ii)}, we have $ \delta_m\ge\gamma>1. $ Since the function $ t\longmapsto\frac{4t}{(1+t)^2} $ is decreasing on $[1,\infty)$, it follows that $$ \mu < \mu_{m,\gamma} \le \mu_{1,\gamma} = \frac{4\gamma}{(1+\gamma)^2}\LN. $$ Hence, Theorem~\ref{thm:main3} applies. 

It remains to prove the additional regularity. Set $ v_n:=u_n^{\frac{1+\delta_m}{2}}. $ If $\delta_m<1$, Lemma~\ref{prioriestimates1} gives $$ \|v_n\|_{H^1_0(\Omega)}\le C, \qquad \|u_n\|_{\sigma_m}\le C. $$ If $\delta_m\ge1$, the same estimates follow from Lemma~\ref{prioriestimates4}. Thus, in either case, $(v_n)$ is bounded in $H^1_0(\Omega)$ and $ \|u_n\|_{\sigma_m}\le C. $ By the weak compactness of $H^1_0(\Omega)$, every subsequence of $(v_n)$ admits a further subsequence and a function $v\in H^1_0(\Omega)$ such that $ v_n\rightharpoonup v $ weakly in $H^1_0(\Omega). $ By the Rellich--Kondrachov theorem, after passing to a further subsequence, $ v_n\to v$ strongly in $L^2(\Omega) $ and a.e. in $\Omega$. On the other hand, since $u_n\to u$ a.e. in $\Omega$, $ v_n = u_n^{\frac{1+\delta_m}{2}} \longrightarrow u^{\frac{1+\delta_m}{2}} $ a.e. in $\Omega. $ Therefore, $ v=u^{\frac{1+\delta_m}{2}} $ a.e. in $\Omega, $ and hence $ u^{\frac{1+\delta_m}{2}}\in H^1_0(\Omega). $ Finally, since $u_n\to u$ a.e. in $\Omega$, Fatou's lemma gives $$ \int_\Omega u^{\sigma_m}\dd x \le \liminf_{n\to\infty} \int_\Omega u_n^{\sigma_m}\dd x \le C. $$ Consequently, $ u\in L^{\sigma_m}(\Omega). $ This completes the proof. \end{proof}

\subsection{Proof of Theorem \ref{thm:uniqueness}}

\begin{lemma}\label{lem:extension}
Let $\gamma\ge0$, $g\in L^1(\Omega)$ nonnegative, $0<\mu<\LN$, and let $u\in H^1_0(\Omega)$ be a weak solution of \eqref{problem1}. Then $g\psi u^{-\gamma}\in L^1(\Omega)$ and \eqref{eq:weakform} hold for every $\psi\in H^1_0(\Omega)$.
\end{lemma}

\begin{proof} For fixed $u\in H^1_0(\Omega)$, the map $ \psi\longmapsto\cE(u,\psi) $ is a bounded linear functional on $H^1_0(\Omega)$. Indeed, by the Cauchy--Schwarz inequality and \eqref{eq:normequiv}, we have  $$ |\cE(u,\psi)| \le \rho(u)\rho(\psi) \le C\|u\|_{H^1_0(\Omega)} \|\psi\|_{H^1_0(\Omega)}. $$ Similarly, by the Cauchy--Schwarz inequality and the Hardy inequality \eqref{eq:hardy}, $$ \begin{aligned} \left| \int_\Omega\frac{u\psi}{|x|^2}\dd x \right| \le \left( \int_\Omega\frac{u^2}{|x|^2}\dd x \right)^{1/2} \left( \int_\Omega\frac{\psi^2}{|x|^2}\dd x \right)^{1/2} \le \frac{1}{\LN} \|\nabla u\|_2\|\nabla\psi\|_2. \end{aligned} $$

Therefore, both functionals are weakly continuous on $H^1_0(\Omega)$. We first extend the weak formulation to bounded test functions with compact support. Let $ \psi\in H^1_0(\Omega)\cap L^\infty(\Omega) $ and suppose that $\supp\psi\subset K$ for some compact set $K\subset\Omega$. Choose a compact neighborhood $K'\subset\Omega$ of $K$. Extending $\psi$ by zero outside $\Omega$, let $ \psi_\varepsilon:=\varrho_\varepsilon*\psi, $ where $(\varrho_\varepsilon)$ is a standard family of mollifiers. For all sufficiently small $\varepsilon>0$, we have $$ \psi_\varepsilon\in C^\infty_c(\Omega), \qquad \supp\psi_\varepsilon\subset K', \qquad \|\psi_\varepsilon\|_\infty\le\|\psi\|_\infty. $$ 

Moreover, $ \psi_\varepsilon\to\psi$ in $H^1_0(\Omega) $ and, along a sequence, $ \psi_\varepsilon\to\psi $ a.e. in $\Omega. $ Since \eqref{eq:weakform} holds for $\psi_\varepsilon$, we have $$ \cE(u,\psi_\varepsilon) - \mu\int_\Omega\frac{u\psi_\varepsilon}{|x|^2}\dd x = \int_\Omega\frac{g\psi_\varepsilon}{u^\gamma}\dd x. $$ The left-hand side converges to $ \cE(u,\psi) - \mu\int_\Omega\frac{u\psi}{|x|^2}\dd x. $ On the other hand, condition \textup{(b)} of Definition~\ref{def:weak} gives a constant $c_{K'}>0$ such that $ u\ge c_{K'}$ a.e. in $K'. $
Consequently, $$ \left| \frac{g\psi_\varepsilon}{u^\gamma} \right| \le c_{K'}^{-\gamma}\|\psi\|_\infty g \quad\text{a.e. in }\Omega. $$ 

The right-hand side belongs to $L^1(\Omega)$. The dominated convergence theorem therefore gives $$ \int_\Omega\frac{g\psi_\varepsilon}{u^\gamma}\dd x \longrightarrow \int_\Omega\frac{g\psi}{u^\gamma}\dd x. $$ Thus, \eqref{eq:weakform} holds for every bounded $\psi\in H^1_0(\Omega)$ with compact support in $\Omega$.

Let now $\psi\in H^1_0(\Omega)$ be nonnegative. Choose $\varphi_k\in C^\infty_c(\Omega)$ such that $ \varphi_k\to\psi $ in $H^1_0(\Omega), $ and define $ \psi_k:=\min\{\psi,\varphi_k^+\}. $ Then $ \psi_k = \psi-(\psi-\varphi_k^+)^+. $ Since the positive-part map is continuous on $H^1_0(\Omega)$ and $\psi^+=\psi$, we have $ \varphi_k^+\to\psi$ in  $H^1_0(\Omega). $ It follows that $ (\psi-\varphi_k^+)^+\to0 $ in $H^1_0(\Omega), $ and hence $ \psi_k\to\psi$ in $H^1_0(\Omega).$  Furthermore, $ 0\le\psi_k\le\psi, $ and $\psi_k$ is bounded and has compact support in $\Omega$, because $ 0\le\psi_k\le\varphi_k^+. $ Thus, the weak formulation already established above applies to $\psi_k$: $$ \int_\Omega\frac{g\psi_k}{u^\gamma}\dd x = \cE(u,\psi_k) - \mu\int_\Omega\frac{u\psi_k}{|x|^2}\dd x. $$ 

Since $\psi_k\to\psi$ strongly in $H^1_0(\Omega)$, the right-hand side converges to $$ \cE(u,\psi) - \mu\int_\Omega\frac{u\psi}{|x|^2}\dd x. $$ After passing to a subsequence, we may also assume that $ \psi_k\to\psi$ a.e. in $\Omega. $ Fatou's lemma then gives $$ \begin{aligned} \int_\Omega\frac{g\psi}{u^\gamma}\dd x &\le \liminf_{k\to\infty} \int_\Omega\frac{g\psi_k}{u^\gamma}\dd x = \cE(u,\psi) - \mu\int_\Omega\frac{u\psi}{|x|^2}\dd x <\infty. \end{aligned} $$ 

Therefore, $ g\psi u^{-\gamma}\in L^1(\Omega). $ Since $$ 0\le\frac{g\psi_k}{u^\gamma} \le \frac{g\psi}{u^\gamma} $$ and $\psi_k\to\psi$ a.e. in $\Omega$, the dominated convergence theorem yields $$ \int_\Omega\frac{g\psi_k}{u^\gamma}\dd x \longrightarrow \int_\Omega\frac{g\psi}{u^\gamma}\dd x. $$ Passing to the limit in the preceding identity, we obtain $$ \cE(u,\psi) - \mu\int_\Omega\frac{u\psi}{|x|^2}\dd x = \int_\Omega\frac{g\psi}{u^\gamma}\dd x. $$ Thus, \eqref{eq:weakform} holds for every nonnegative $\psi\in H^1_0(\Omega)$.

Finally, let $\psi\in H^1_0(\Omega)$ be arbitrary. Write $$ \psi=\psi^+-\psi^-, \qquad \psi^+:=\max\{\psi,0\}, \qquad \psi^-:=\max\{-\psi,0\}. $$ Since $\psi^\pm\in H^1_0(\Omega)$, the preceding argument gives $ g\psi^\pm u^{-\gamma}\in L^1(\Omega) $ and \eqref{eq:weakform} for both $\psi^+$ and $\psi^-$. Consequently, $ g\psi u^{-\gamma}\in L^1(\Omega), $ and subtracting the two corresponding identities proves \eqref{eq:weakform} for $\psi$. \end{proof}

\begin{proof}[Proof of Theorem~\ref{thm:uniqueness}] The argument follows the one used for finite-energy solutions of singular problems involving the Laplacian; see \cite{boccardo2014some,canino2004variational}. The first assertion follows from Lemma~\ref{lem:extension}. Let $u,v\in H^1_0(\Omega)$ be two weak solutions of \eqref{problem1}, and set $ z:=u-v. $ By Lemma~\ref{lem:extension}, the weak formulation is valid for every test function in $H^1_0(\Omega)$. We may therefore use $z^+\in H^1_0(\Omega)$ as a test function in the equations satisfied by $u$ and $v$. Subtracting the resulting identities, we obtain $$ \begin{aligned} \cE(z,z^+) -\mu\int_\Omega\frac{zz^+}{|x|^2}\dd x &= \int_\Omega g\left( \frac{1}{u^\gamma} -\frac{1}{v^\gamma} \right)z^+\dd x \\ &\le0. \end{aligned} $$ 

Indeed, on the set $\{z>0\}=\{u>v\}$, the function $t\mapsto t^{-\gamma}$ is nonincreasing, and hence $ \frac{1}{u^\gamma}-\frac{1}{v^\gamma}\le0. $ For $\gamma=0$, the difference on the left-hand side of this inequality is identically zero. On the other hand, Lemma~\ref{lem:alg}\textup{(iii)} gives $ \cE(z,z^+)\ge\rho(z^+)^2. $ Since $ zz^+=(z^+)^2 $ a.e. in $\Omega, $ the Hardy inequality \eqref{eq:hardy} yields $$ \begin{aligned} \cE(z,z^+) -\mu\int_\Omega\frac{zz^+}{|x|^2}\dd x &\ge \rho(z^+)^2 -\mu\int_\Omega\frac{(z^+)^2}{|x|^2}\dd x \\ &\ge \rho(z^+)^2 -\frac{\mu}{\LN}\|\nabla z^+\|_2^2 \ge \left(1-\frac{\mu}{\LN}\right)\rho(z^+)^2. \end{aligned} $$ 

Combining the preceding inequalities, we obtain $ \left(1-\frac{\mu}{\LN}\right)\rho(z^+)^2\le0. $ Since $\mu<\LN$, it follows that $ \rho(z^+)=0. $ Thus $z^+=0$ a.e. in $\Omega$, and consequently $ u\le v $ a.e. in $\Omega. $ Exchanging the roles of $u$ and $v$, we obtain $ v\le u$  a.e. in $\Omega. $ Therefore, $ u=v$ a.e. in $\Omega. $ Finally, the solution obtained in Theorem~\ref{thm:main2} belongs to $H^1_0(\Omega)$. Hence it is the unique finite-energy weak solution of \eqref{problem1}, which proves the last assertion. \end{proof}

\subsection{Proof of Theorem \ref{thm:nonexistence}}
We shall use Picone-type inequalities for both the local and nonlocal parts of $\cL$. Part \textup{(b)} is the classical Picone identity; see \cite{allegretto1998picone} and the references therein. Part \textup{(a)} is the elementary algebraic inequality underlying the nonlocal Picone inequality and the ground-state representation for the fractional Laplacian; see \cite{brasco2014convexity,frank2008nonlinear}. We include the proofs in the required form below, ensuring that the function $u$ vanishes outside $\Omega$.

\begin{lemma}\label{lem:picone}
\begin{enumerate}[label=(\alph*)]
\item Let $a,b\ge0$ and $c,d\in\R$, and assume that $a>0$ if $c\ne0$, and $b>0$ if $d\neq0$. Set $c^2/a:=0$ if $c=0$ and $d^2/b:=0$ if $d=0$. Then
\[
(a-b)\Big(\frac{c^2}{a}-\frac{d^2}{b}\Big)\le(c-d)^2 .
\]

\item Let $U\subset\RN$ be open, $u\in H^1(U)$ with $u\ge c_0>0$ a.e. in $U$, and $\eta\in C^\infty_c(U)$. Then $\eta^2/u\in H^1(\RN)\cap L^\infty(\RN)$ has compact support in $U$, and
\[
\nabla u\cdot\nabla\Big(\frac{\eta^2}{u}\Big)=|\nabla\eta|^2-\Big|\nabla\eta-\frac{\eta}{u}\nabla u\Big|^2\le|\nabla\eta|^2\qquad\text{a.e. in }U.
\]
\end{enumerate}
\end{lemma}

\begin{proof}
(a) If $c\ne0$ and $d\neq0$, then $a,b>0$ and
\[
(a-b)\Big(\frac{c^2}{a}-\frac{d^2}{b}\Big)=c^2+d^2-\Big(\frac ba\,c^2+\frac ab\,d^2\Big)\le c^2+d^2-2|c||d|\le(c-d)^2,
\]
by the arithmetic--geometric mean inequality. If $d=0$, the left-hand side equals $c^2-b\,c^2/a\le c^2$ (or $0$ if also $c=0$), since $b\ge0$; the case $c=0$ is symmetric.

(b) Since $u\ge c_0>0$ a.e. in $U$, we have $ \frac1u = \frac{1}{\max\{u,c_0\}}$ a.e. in $U$. The function $ t\longmapsto\frac{1}{\max\{t,c_0\}} $ is bounded and Lipschitz on $\R$. The  chain rule therefore gives $ \frac1u\in H^1_{\rm loc}(U)\cap L^\infty(U) $ and $$ \nabla\left(\frac1u\right) = -\frac{\nabla u}{u^2} \quad\text{a.e. in }U. $$ 

Since $\eta$ has compact support in $U$, it follows that $ \frac{\eta^2}{u}\in H^1_0(U)\cap L^\infty(U). $ Its zero extension outside $U$ consequently belongs to $H^1(\RN)\cap L^\infty(\RN)$ and has compact support in $U$. By the product and chain rules, $ \nabla\left(\frac{\eta^2}{u}\right) = \frac{2\eta}{u}\nabla\eta - \frac{\eta^2}{u^2}\nabla u $ a.e. in $U$. Hence, $$ \nabla u\cdot\nabla\left(\frac{\eta^2}{u}\right) = \frac{2\eta}{u}\nabla u\cdot\nabla\eta - \frac{\eta^2}{u^2}|\nabla u|^2. $$ On the other hand, $$ \left| \nabla\eta-\frac{\eta}{u}\nabla u \right|^2 = |\nabla\eta|^2 - \frac{2\eta}{u}\nabla u\cdot\nabla\eta + \frac{\eta^2}{u^2}|\nabla u|^2. $$ Combining the last two identities yields $$ \nabla u\cdot\nabla\left(\frac{\eta^2}{u}\right) = |\nabla\eta|^2 - \left| \nabla\eta-\frac{\eta}{u}\nabla u \right|^2 \le |\nabla\eta|^2 $$ a.e. in $U$.
\end{proof}

\begin{proof}[Proof of Theorem~\ref{thm:nonexistence}] Assume by contradiction that there exists a function $u$ satisfying the assumptions of the theorem. Since $N\ge3$, we have $ \frac{N}{N-1}<2<N, $ so $q=2$ is admissible in Lemma~\ref{lem:welldef}. Consequently, all the terms in \eqref{eq:supersol} are well defined. We first extend \eqref{eq:supersol} to nonnegative $H^1$ test functions with compact support in $\Omega$. Fix a compact set $K'\subset\Omega$. Choose open sets $ K'\Subset U'\Subset U\Subset\Omega $ and a function $\varphi\in C^\infty_c(U)$ such that $ 0\le\varphi\le1$ and $ \varphi\equiv1 $ in $U'. $ 

Set $ u_1:=u\varphi, \, u_2:=u(1-\varphi). $ Then $u_1\in H^1_0(U)$, while $u_2=0$ in $U'$. Set $ d:={\rm dist}(K',\RN\setminus U')>0. $ As in the proof of Lemma~\ref{lem:welldef}, for every $\psi\in H^1(\RN)$ with $\supp\psi\subset K'$, we have $$ \begin{aligned} \cE(u,\psi) &= \int_U\nabla u\cdot\nabla\psi\dd x + \frac{\Cs}{2} \iint_{\R^{2N}} \frac{ (u_1(x)-u_1(y))(\psi(x)-\psi(y)) } {|x-y|^{N+2s}}\dd x\dd y \\ &\quad\quad\quad - \Cs\int_{K'}\psi(x) \int_{\RN\setminus U'} \frac{u_2(y)}{|x-y|^{N+2s}}\dd y\dd x. \end{aligned} $$ 

All the integrals in this identity converge absolutely. Indeed, the absolute value of the second term is bounded by $ \frac{\Cs}{2}[u_1]_s[\psi]_s, $ whereas the absolute value of the last term is bounded by $ \Cs d^{-N-2s} \|u\|_{L^1(\Omega)} \|\psi\|_{L^1(\Omega)}. $ Since $0<s<1$, the embedding $ H^1(\RN)\hookrightarrow H^s(\RN) $ is continuous. It follows that the map $ \psi\longmapsto\cE(u,\psi) $ is continuous with respect to $H^1(\RN)$ convergence among functions supported in $K'$. The map $$ \psi\longmapsto \int_\Omega\frac{u\psi}{|x|^2}\dd x $$ is also continuous on this class of functions. Indeed, since $\varphi=1$ on $K'$, we have $u\psi=u_1\psi$, and therefore, by the Cauchy--Schwarz and Hardy inequalities, $$ \begin{aligned} \left| \int_\Omega\frac{u\psi}{|x|^2}\dd x \right| = \left| \int_\Omega\frac{u_1\psi}{|x|^2}\dd x \right| \le \left( \int_\Omega\frac{u_1^2}{|x|^2}\dd x \right)^{1/2} \left( \int_\Omega\frac{\psi^2}{|x|^2}\dd x \right)^{1/2} \le \frac{1}{\LN} \|\nabla u_1\|_2\|\nabla\psi\|_2. \end{aligned} $$ 

Now let $\psi\in H^1(\RN)$ be nonnegative and have compact support $K$ in the interior of $K'$. Extending $\psi$ by zero where necessary, define $ \psi_\varepsilon:=\varrho_\varepsilon*\psi, $ where $(\varrho_\varepsilon)$ is a standard family of nonnegative mollifiers. For all sufficiently small $\varepsilon>0$, we have $$ \psi_\varepsilon\in C^\infty_c(\Omega), \qquad \psi_\varepsilon\ge0, \qquad \supp\psi_\varepsilon\subset K', $$ and $ \psi_\varepsilon\to\psi$ in $H^1(\RN). $ Applying \eqref{eq:supersol} to $\psi_\varepsilon$ and passing to the limit, we conclude that $$ \cE(u,\psi) - \mu\int_\Omega\frac{u\psi}{|x|^2}\dd x \ge0 $$ for every nonnegative $\psi\in H^1(\RN)$ with compact support in $\Omega$. Let now $\eta\in C^\infty_c(\Omega)$ and set $ K:=\supp\eta. $ Choose an open set $U$ such that $ K\Subset U\Subset\Omega. $ By condition \textup{(b)} of Definition~\ref{def:weak}, there exists $c_U>0$ such that $ u\ge c_U$ a.e. in $U. $ Lemma~\ref{lem:picone}\textup{(b)} shows that the zero extension of $ \psi:=\frac{\eta^2}{u} $ belongs to $H^1(\RN)\cap L^\infty(\RN)$, is nonnegative, and has compact support in $U$. Hence it is an admissible test function in the extended supersolution inequality. Moreover, $$ \nabla u\cdot\nabla\left(\frac{\eta^2}{u}\right) \le |\nabla\eta|^2 \quad\text{a.e. in }U. $$ For the nonlocal term, Lemma~\ref{lem:picone}\textup{(a)}, applied with $$ a=u(x), \qquad b=u(y), \qquad c=\eta(x), \qquad d=\eta(y), $$ gives $$ (u(x)-u(y)) \left( \frac{\eta(x)^2}{u(x)} - \frac{\eta(y)^2}{u(y)} \right) \le (\eta(x)-\eta(y))^2 $$ for a.e. $x,y\in\RN$. The hypotheses of Lemma~\ref{lem:picone}\textup{(a)} are satisfied because $u\ge0$ in $\RN$ and $u>0$ wherever $\eta\ne0$. We use the convention $ \frac{\eta(x)^2}{u(x)}:=0 \quad\text{whenever }\eta(x)=0. $ 

It follows that $$ \cE\left(u,\frac{\eta^2}{u}\right) \le \int_{\RN}|\nabla\eta|^2\dd x + \frac{\Cs}{2}[\eta]_s^2 = \rho(\eta)^2. $$ Since $ u\frac{\eta^2}{u}=\eta^2$ a.e. in $\Omega, $ the extended supersolution inequality gives $$ \mu\int_\Omega\frac{\eta^2}{|x|^2}\dd x \le \cE\left(u,\frac{\eta^2}{u}\right) \le \rho(\eta)^2. $$ Therefore, $$ \mu \le \frac{\rho(\eta)^2} {\displaystyle\int_\Omega\frac{\eta^2}{|x|^2}\dd x} $$ for every nonzero $\eta\in C^\infty_c(\Omega)$. Taking the infimum over all such functions and using Remark~\ref{rem:hardy-optimal}, we obtain $$ \mu \le \inf_{\eta\in C^\infty_c(\Omega)\setminus\{0\}} \frac{\rho(\eta)^2} {\displaystyle\int_\Omega\frac{\eta^2}{|x|^2}\dd x} = \LN. $$ 

This contradicts the assumption $\mu>\LN$. Finally, if $u$ were a weak solution of \eqref{problem1}, then, for every nonnegative $\psi\in C^\infty_c(\Omega)$, $$ \begin{aligned} \cE(u,\psi) - \mu\int_\Omega\frac{u\psi}{|x|^2}\dd x = \int_\Omega\frac{g\psi}{u^\gamma}\dd x \ge0, \end{aligned} $$ because $g\ge0$. Thus $u$ would satisfy \eqref{eq:supersol}, which has just been shown to be impossible. Hence \eqref{problem1} admits no weak solution. \end{proof}

\begin{remark}\label{rem:nonexistence} The proof of Theorem~\ref{thm:nonexistence} yields the following additional information. Let $\mu>0$, and let $u\in H^1_{\rm loc}(\Omega)$ be a weak solution of \eqref{problem1}. Then, for every $\eta\in C^\infty_c(\Omega)$, the weak formulation \eqref{eq:weakform} remains valid with the test function $ \psi:=\frac{\eta^2}{u}. $ Indeed, condition \textup{(b)} of Definition~\ref{def:weak} implies that $u$ is bounded away from zero on $\supp\eta$. Hence $\psi$ is bounded, nonnegative, and compactly supported in $\Omega$. The extension follows by the same approximation argument used in the proof of Theorem~\ref{thm:nonexistence}. In particular, the singular term passes to the limit by the dominated convergence theorem. Using $\psi=\eta^2/u$ in \eqref{eq:weakform} gives $$ \cE\left(u,\frac{\eta^2}{u}\right) - \mu\int_\Omega\frac{\eta^2}{|x|^2}\dd x = \int_\Omega\frac{g\eta^2}{u^{1+\gamma}}\dd x. $$ By Lemma~\ref{lem:picone}, $$ \cE\left(u,\frac{\eta^2}{u}\right) \le \rho(\eta)^2. $$ Consequently, $$ \mu\int_\Omega\frac{\eta^2}{|x|^2}\dd x + \int_\Omega\frac{g\eta^2}{u^{1+\gamma}}\dd x \le \rho(\eta)^2 $$ for every $\eta\in C^\infty_c(\Omega)$. In particular, both integrals on the left-hand side are finite. This weighted Hardy-type inequality may provide a useful starting point for the analysis of the critical case $\mu=\LN$. \end{remark}

\section{Behaviour near the origin and the critical case}\label{sec:origin}

Throughout this section we write $f(t):=t(N-2-t)$, so that $-\Delta|x|^{-t}=f(t)|x|^{-t-2}$ in $\RN\setminus\{0\}$ for every $t\in\R$; $f$ is increasing on $(-\infty,\frac{N-2}{2}]$, and $\alpha_\mu$ defined in \eqref{eq:alphamu} is the unique solution of $f(t)=\mu$ in $(0,\frac{N-2}{2}]$.

\subsection{An explicit subsolution}
In the local case, the lower bound $u\ge A|x|^{-\alpha_\mu}$ is obtained by comparison with explicit radial subsolutions; see \cite{brezis2005semilinear, dupaigne2002nonlinear}. For $\cL$, the nonlocal part of such a subsolution has to be controlled; to this end, we first record two elementary facts on the fractional Laplacian. The exact value of the constant $c_\beta$ in Lemma \ref{lem:fraclap}(a) is computed in \cite{frank2008hardy}, but we only need its finiteness.

\begin{lemma}\label{lem:fraclap}
\begin{enumerate}[label=(\alph*)]
\item Let $0<\beta<N-2s$ and $h_\beta(x):=|x|^{-\beta}$. For every $x\neq0$, the principal value
\[
(-\Delta)^sh_\beta(x)=\Cs\lim_{\varepsilon\to0^+}\int_{\{|x-y|\ge\varepsilon\}}\frac{h_\beta(x)-h_\beta(y)}{|x-y|^{N+2s}}\dd y
\]
exists, and $(-\Delta)^sh_\beta(x)=c_\beta|x|^{-\beta-2s}$ for some constant $c_\beta\in\R$ depending only on $N,s,\beta$.
\item Let $Z\in L^1(\RN)\cap H^s(\RN)$, and assume that $Z\in C^2(V)$ for some open set $V\subset\RN$. Then the principal value $(-\Delta)^sZ(x)$ exists for every $x\in V$, it is locally bounded in $V$, and
\[
\frac{\Cs}{2}\iint_{\R^{2N}}\frac{(Z(x)-Z(y))(\psi(x)-\psi(y))}{|x-y|^{N+2s}}\dd x\dd y=\int_{\RN}(-\Delta)^sZ(x)\,\psi(x)\dd x
\]
for all $\psi\in C^\infty_c(V)$.
\end{enumerate}
\end{lemma}

\begin{proof}

\textup{(a)} Fix $x\ne0$ and choose $r\in(0,|x|/2)$. The function $h_\beta$ is smooth in $B_r(x)$. By the change of variables $y=x+z$ and the symmetry of the annulus, we have $$ \begin{aligned} &\int_{\{\varepsilon\le|x-y|<r\}} \frac{h_\beta(x)-h_\beta(y)} {|x-y|^{N+2s}}\dd y \\ &\qquad= -\frac12 \int_{\{\varepsilon\le|z|<r\}} \frac{ h_\beta(x+z)+h_\beta(x-z)-2h_\beta(x) } {|z|^{N+2s}}\dd z. \end{aligned} $$ By Taylor's theorem, $ \left| h_\beta(x+z)+h_\beta(x-z)-2h_\beta(x) \right| \le C_{x,r}|z|^2 $ for $|z|<r$, where $ C_{x,r} := \|D^2h_\beta\|_{L^\infty(B_r(x))}. $ Since $s<1$, $$ \int_{\{|z|<r\}} \frac{|z|^2}{|z|^{N+2s}}\dd z<\infty. $$ Therefore, the integral over $\{|x-y|<r\}$ converges as $\varepsilon\to0^+$. The integral over $\{|x-y|\ge r\}$ converges absolutely. Indeed, $h_\beta\in L^1_{\rm loc}(\RN)$ because $\beta<N$, while at infinity $$ \frac{|y|^{-\beta}}{|x-y|^{N+2s}} \le C|y|^{-N-2s-\beta}, $$ which is integrable because $\beta+2s>0$. Thus, the principal value $(-\Delta)^sh_\beta(x)$ exists for every $x\ne0$. To determine its homogeneity, write $ x=\lambda\xi, \, \lambda:=|x|, \, \xi:=\frac{x}{|x|}, $ and make the change of variables $y=\lambda z$. This gives $$ (-\Delta)^sh_\beta(x) = |x|^{-\beta-2s} (-\Delta)^sh_\beta\left(\frac{x}{|x|}\right). $$ Since both $h_\beta$ and the kernel are invariant under rotations, the value of $(-\Delta)^sh_\beta$ is constant on the unit sphere. Denoting this value by $c_\beta$, we obtain $ (-\Delta)^sh_\beta(x) = c_\beta|x|^{-\beta-2s}. $

\textup{(b)} Fix a compact set $K\Subset V$ and choose $r>0$ such that $ \{x\in\RN:{\rm dist}(x,K)\le2r\}\Subset V. $ For $x\in K$, split the principal-value integral into its contributions over $B_r(x)$ and $\RN\setminus B_r(x)$. By symmetrization, $$ \begin{aligned} \int_{\{\varepsilon\le|x-y|<r\}} \frac{Z(x)-Z(y)} {|x-y|^{N+2s}}\dd y = -\frac12 \int_{\{\varepsilon\le|z|<r\}} \frac{Z(x+z)+Z(x-z)-2Z(x)} {|z|^{N+2s}}\dd z. \end{aligned} $$ Since $Z\in C^2(V)$, Taylor's theorem gives $ |Z(x+z)+Z(x-z)-2Z(x)| \le C_K|z|^2 $ for $x\in K$ and $|z|<r$, where $C_K$ is independent of $x$. Since $s<1$, the resulting majorant is integrable near the origin. For the contribution over $\RN\setminus B_r(x)$, we have $$ \begin{aligned} \int_{\RN\setminus B_r(x)} \frac{|Z(x)-Z(y)|} {|x-y|^{N+2s}}\dd y \le |Z(x)| \int_{\{|z|\ge r\}} \frac{\dd z}{|z|^{N+2s}} + r^{-N-2s}\|Z\|_{L^1(\RN)}. \end{aligned} $$

Since $Z$ is bounded on $K$, this estimate is uniform for $x\in K$. It follows that the principal value exists locally uniformly in $V$ and that $(-\Delta)^sZ$ is locally bounded there. Let now $\psi\in C^\infty_c(V)$ and $\delta>0$. Since the kernel is bounded on $\{|x-y|\ge\delta\}$ and $Z\in L^1(\RN)$, all the following integrals converge absolutely. By expanding the product, exchanging $x$ and $y$ in one of the resulting terms, and using the symmetry of the kernel, we obtain $$ \begin{aligned} &\frac{\Cs}{2} \iint_{\{|x-y|\ge\delta\}} \frac{(Z(x)-Z(y))(\psi(x)-\psi(y))} {|x-y|^{N+2s}}\dd x\dd y \\ &\qquad= \Cs\int_{\RN}\psi(x) \int_{\{|x-y|\ge\delta\}} \frac{Z(x)-Z(y)} {|x-y|^{N+2s}}\dd y\dd x. \end{aligned} $$

As $\delta\to0^+$, the integrand on the left is dominated by $$ \frac{ |Z(x)-Z(y)||\psi(x)-\psi(y)| } {|x-y|^{N+2s}}, $$ which belongs to $L^1(\R^{2N})$ by the Cauchy--Schwarz inequality, because $Z,\psi\in H^s(\RN)$. Hence, the dominated convergence theorem gives $$ \begin{aligned} &\lim_{\delta\to0^+} \frac{\Cs}{2} \iint_{\{|x-y|\ge\delta\}} \frac{(Z(x)-Z(y))(\psi(x)-\psi(y))} {|x-y|^{N+2s}}\dd x\dd y \\ &\qquad= \frac{\Cs}{2} \iint_{\R^{2N}} \frac{(Z(x)-Z(y))(\psi(x)-\psi(y))} {|x-y|^{N+2s}}\dd x\dd y. \end{aligned} $$ 

On the right-hand side, the truncated inner integral converges uniformly for $x\in\supp\psi$ to $ \Cs^{-1}(-\Delta)^sZ(x). $ Since $(-\Delta)^sZ$ is bounded on $\supp\psi$, we may pass to the limit under the outer integral. We conclude that $$ \frac{\Cs}{2} \iint_{\R^{2N}} \frac{(Z(x)-Z(y))(\psi(x)-\psi(y))} {|x-y|^{N+2s}}\dd x\dd y = \int_{\RN}(-\Delta)^sZ(x)\psi(x)\dd x, $$
concluding the proof.
\end{proof}

\begin{lemma}\label{lem:subsol}
Let $0<\mu<\LN$, $\alpha:=\alpha_\mu$, and fix $\tau\in(0,\min\{\alpha,2-2s\})$. Let
\[
\Theta(x):=|x|^{-\alpha}+|x|^{-\alpha+\tau},\qquad Z_R:=\Theta\,\varphi_R,
\]
where $R\in(0,\frac13]$ and $\varphi_R\in C^\infty_c(B_{3R})$ satisfies $0\le\varphi_R\le1$ and $\varphi_R\equiv1$ in $B_{2R}$. Then $Z_R\in H^1(\RN)\cap L^1(\RN)$, and there exists $R_0\in(0,\frac13]$, depending only on $N,s,\mu,\tau$, such that for every $R\in(0,R_0]$
\begin{equation}\label{eq:subsol}
\cE(Z_R,\psi)-\mu\int_{B_R}\frac{Z_R\,\psi}{|x|^2}\dd x\le0\qquad\text{for every }\psi\in H^1_0(B_R),\ \psi\ge0 .
\end{equation}
\end{lemma}
\begin{proof} Since $0<\mu<\LN$ and $ 0<\tau<\alpha<\frac{N-2}{2}, $ we have $ 0<\alpha-\tau<\alpha<\frac{N-2}{2}. $ In particular, $ |\Theta(x)|\le C|x|^{-\alpha} $ and $ |\nabla\Theta(x)|\le C|x|^{-\alpha-1} $ near the origin. Since $ 2\alpha+2<N, $ it follows that $ \Theta\in H^1_{\rm loc}(\RN)\cap L^1_{\rm loc}(\RN). $ As $\varphi_R$ is smooth and compactly supported, we conclude that $ Z_R\in H^1(\RN)\cap L^1(\RN). $ Moreover, $ Z_R\in C^\infty(\RN\setminus\{0\}). $

Recall that $$ -\Delta|x|^{-\beta} = \beta(N-2-\beta)|x|^{-\beta-2} = f(\beta)|x|^{-\beta-2} $$ for $0<\beta<N-2$. Since $f(\alpha)=\mu$, we obtain that, in $B_{2R}\setminus\{0\},$ the following holds $$ \begin{aligned} -\Delta\Theta-\mu\frac{\Theta}{|x|^2} = \bigl(f(\alpha)-\mu\bigr)|x|^{-\alpha-2} + \bigl(f(\alpha-\tau)-\mu\bigr) |x|^{-\alpha+\tau-2} = -\kappa|x|^{-\alpha+\tau-2}, \end{aligned} $$ where $ \kappa:=\mu-f(\alpha-\tau)>0. $ Indeed, $f$ is strictly increasing on $(0,(N-2)/2]$, and $ 0<\alpha-\tau<\alpha<\frac{N-2}{2}. $ Since $\alpha<N-2$, the preceding identity also holds in the sense of distributions in $B_{2R}$. We next estimate the fractional part. Let $ \beta\in\{\alpha,\alpha-\tau\}. $ Then $ 0<\beta<\frac{N-2}{2}<N-2s. $ For $0<|x|<R$, write $$ h_\beta\varphi_R=h_\beta-w_R, \qquad w_R:=h_\beta(1-\varphi_R)\ge0. $$ Since $w_R=0$ in $B_{2R}$, no principal value is needed when computing $(-\Delta)^sw_R(x)$ for $|x|<R$. By Lemma~\ref{lem:fraclap}\textup{(a)}, we have $$ \begin{aligned} (-\Delta)^s(h_\beta\varphi_R)(x) = c_\beta|x|^{-\beta-2s}+ \Cs\int_{\{|y|\ge2R\}} \frac{w_R(y)}{|x-y|^{N+2s}}\dd y. \end{aligned} $$ 

For $|x|<R$ and $|y|\ge2R$, we have $ |x-y|\ge|y|-|x|\ge\frac{|y|}{2}. $ Since $ 0\le w_R(y)\le|y|^{-\beta}, $ it follows that $$ \begin{aligned} \int_{\{|y|\ge2R\}} \frac{w_R(y)}{|x-y|^{N+2s}}\dd y \le 2^{N+2s} \int_{\{|y|\ge2R\}} |y|^{-N-2s-\beta}\dd y \le CR^{-\beta-2s} \le C|x|^{-\beta-2s}. \end{aligned} $$ 

Therefore, $ (-\Delta)^s(h_\beta\varphi_R)(x) \le C_\beta|x|^{-\beta-2s} $ for $0<|x|<R$, where $C_\beta$ is independent of $R$. Since $R\le1/3$ and $ |x|^{-\alpha+\tau-2s} \le |x|^{-\alpha-2s} $ for $0<|x|<R$, we conclude that \begin{equation}\label{eq:nonlocal-up} (-\Delta)^sZ_R(x) \le \widehat C|x|^{-\alpha-2s} \end{equation} for $0<|x|<R$, where $\widehat C>0$ depends only on $N$, $s$, $\alpha$, and $\tau$, and is independent of $R$.

Let now $ \psi\in C^\infty_c(B_R\setminus\{0\}), \, \psi\ge0. $ Since $Z_R=\Theta$ in $B_{2R}$, integration by parts in the local term, together with Lemma~\ref{lem:fraclap}\textup{(b)}, gives $$ \begin{aligned} \cE(Z_R,\psi) -\mu\int_{B_R}\frac{Z_R\psi}{|x|^2}\dd x  = \int_{B_R} \left( -\kappa|x|^{-\alpha+\tau-2} + (-\Delta)^sZ_R \right)\psi\dd x. \end{aligned} $$ 

Using \eqref{eq:nonlocal-up}, we obtain $$ \begin{aligned} \cE(Z_R,\psi) -\mu\int_{B_R}\frac{Z_R\psi}{|x|^2}\dd x\le \int_{B_R} |x|^{-\alpha-2} \left( \widehat C|x|^{2-2s} -\kappa|x|^\tau \right)\psi\dd x. \end{aligned} $$ 

Since $\tau<2-2s$, define $$ R_0 := \min\left\{ \frac13, \left(\frac{\kappa}{\widehat C}\right)^{ \frac{1}{2-2s-\tau} } \right\}. $$ If $0<|x|\le R\le R_0$, then $ \widehat C|x|^{2-2s} \le \kappa|x|^\tau. $ Consequently, we have  $$ \cE(Z_R,\psi) - \mu\int_{B_R}\frac{Z_R\psi}{|x|^2}\dd x \le0 $$ for every nonnegative $\psi\in C^\infty_c(B_R\setminus\{0\})$. It remains to extend this inequality to nonnegative functions in $H^1_0(B_R)$. The map $$ \psi\longmapsto \cE(Z_R,\psi) - \mu\int_{B_R}\frac{Z_R\psi}{|x|^2}\dd x $$ is a continuous linear functional on $H^1_0(B_R)$. Indeed, $Z_R\in H^1(\RN)$, and therefore $ |\cE(Z_R,\psi)| \le \rho(Z_R)\rho(\psi). $ Moreover, $ Z_R/|x| \in L^2(B_R), $ because $2\alpha+2<N$, and the Hardy inequality gives $$ \begin{aligned} \left| \int_{B_R}\frac{Z_R\psi}{|x|^2}\dd x \right| \le \left\|\frac{Z_R}{|x|}\right\|_2 \left\|\frac{\psi}{|x|}\right\|_2 \le \frac{1}{\sqrt{\LN}} \left\|\frac{Z_R}{|x|}\right\|_2 \|\nabla\psi\|_2. \end{aligned} $$ 

Since $N\ge3$, the point $\{0\}$ has zero $H^1$-capacity. Hence $ C^\infty_c(B_R\setminus\{0\}) $ is dense in $H^1_0(B_R)$. More precisely, every nonnegative $\psi\in H^1_0(B_R)$ can be approximated in $H^1_0(B_R)$ by nonnegative functions $ \psi_k\in C^\infty_c(B_R\setminus\{0\}). $ For completeness, let $\zeta\in C^\infty_c(B_2)$ satisfy $ 0\le\zeta\le1$ and $\zeta\equiv1$ in $B_1, $ and define $ \widetilde\psi_k := \psi\bigl(1-\zeta(k\,\cdot)\bigr). $ Then $ \widetilde\psi_k\to\psi $ in $L^2(B_R), $ and $ \zeta(k\,\cdot)\nabla\psi\to0 $ in $L^2(B_R). $
Furthermore, $$ k|\nabla\zeta(kx)| \le \frac{2\|\nabla\zeta\|_\infty}{|x|} $$ on the support of $\nabla\zeta(k\,\cdot)$. Thus, $$ \begin{aligned} \left\| k\psi\nabla\zeta(k\,\cdot) \right\|_2 \le 2\|\nabla\zeta\|_\infty \left\| \frac{\psi}{|x|} \right\|_{L^2(B_{2/k})} \longrightarrow0. \end{aligned} $$ Therefore, $ \widetilde\psi_k\to\psi \quad\text{in }H^1_0(B_R). $ Each $\widetilde\psi_k$ is nonnegative and vanishes in a neighborhood of the origin. Standard nonnegative approximation in the Lipschitz domain $ B_R\setminus\overline{B_{1/k}} $ then yields nonnegative functions in $C^\infty_c(B_R\setminus\{0\})$ converging to $\widetilde\psi_k$ in $H^1$. By continuity of the preceding linear functional, \eqref{eq:subsol} passes to the limit and holds for every nonnegative $\psi\in H^1_0(B_R)$. \end{proof}

\subsection{Proof of Theorem \ref{thm:origin}}
\begin{proof}[Proof of Theorem~\ref{thm:origin}\textup{(a)}] Let $\tau$ and $R_0$ be as in Lemma~\ref{lem:subsol}. Choose $R\in(0,R_0]$ such that $ \overline{B_{4R}}\subset\Omega. $ By condition \textup{(b)} of Definition~\ref{def:weak}, there exists $c>0$ such that $ u\ge c$ a.e. in $B_{3R}.$ Set $ A:=\frac{cR^\alpha}{2}, \, Z:=AZ_R, $ where $\alpha=\alpha_\mu$. If $R\le|x|<3R$, then, since $R\le1/3$, we obtain $$ \begin{aligned} Z(x) = A\varphi_R(x) \left( |x|^{-\alpha}+|x|^{-\alpha+\tau} \right) \le A|x|^{-\alpha}(1+|x|^\tau) \le 2AR^{-\alpha} = c \le u(x). \end{aligned} $$ 

Moreover, $Z=0\le u$ a.e. in $\RN\setminus B_{3R}$. Therefore, $ Z\le u$ a.e. in $\RN\setminus B_R. $ Choose $\eta\in C^\infty_c(B_{4R})$ such that $ 0\le\eta\le1$ and $\eta\equiv1$ in $B_{3R},$ and define $ \widetilde u:=u\eta. $ Since $u\in H^1_{\rm loc}(\Omega)$, we have $ \widetilde u\in H^1(\RN). $ Furthermore, $Z$ vanishes outside $B_{3R}$ and $\widetilde u=u$ in $B_{3R}$. Hence $ (Z-u)^+=(Z-\widetilde u)^+ $ a.e. in $\RN. $ It follows that $ z^+:=(Z-u)^+\in H^1(\RN). $ Since $Z\le u$ a.e. in $\RN\setminus B_R$, we also have $ z^+=0$ a.e. in $\RN\setminus B_R. $ Therefore, $ z^+\in H^1_0(B_R). $ 

The approximation argument used in the proof of Theorem~\ref{thm:nonexistence}, which does not depend on the assumption $\mu>\LN$, extends \eqref{eq:supersol} to the nonnegative test function $z^+$. Thus, $$ \cE(u,z^+) - \mu\int_{B_R}\frac{uz^+}{|x|^2}\dd x \ge0. $$ On the other hand, Lemma~\ref{lem:subsol}, applied to $z^+$, gives $$ \cE(Z,z^+) - \mu\int_{B_R}\frac{Zz^+}{|x|^2}\dd x \le0. $$ Subtracting the first inequality from the second and setting $ z:=Z-u, $ we obtain $$ \begin{aligned} \int_{\RN}\nabla z\cdot\nabla z^+\dd x + \frac{\Cs}{2} \iint_{\R^{2N}} \frac{ (z(x)-z(y))(z^+(x)-z^+(y)) } {|x-y|^{N+2s}}\dd x\dd y - \mu\int_{B_R}\frac{zz^+}{|x|^2}\dd x \le0. \end{aligned} $$ 

All the terms are well defined by the preceding extension argument and Lemma~\ref{lem:subsol}. Since $ \nabla z\cdot\nabla z^+=|\nabla z^+|^2 $ and $ zz^+=(z^+)^2 $ a.e., while Lemma~\ref{lem:alg}\textup{(iii)} gives $$ (z(x)-z(y))(z^+(x)-z^+(y)) \ge (z^+(x)-z^+(y))^2, $$ we conclude that $$ \rho(z^+)^2 - \mu\int_{B_R}\frac{(z^+)^2}{|x|^2}\dd x \le0. $$ 

By the Hardy inequality, $$ \begin{aligned} \rho(z^+)^2 - \mu\int_{B_R}\frac{(z^+)^2}{|x|^2}\dd x \ge \rho(z^+)^2 - \frac{\mu}{\LN}\|\nabla z^+\|_2^2 \ge \left(1-\frac{\mu}{\LN}\right)\rho(z^+)^2. \end{aligned} $$ 

Therefore, $$ \left(1-\frac{\mu}{\LN}\right)\rho(z^+)^2\le0. $$ 

Since $\mu<\LN$, it follows that $ z^+=0. $ Thus, $ u\ge Z = A\left( |x|^{-\alpha}+|x|^{-\alpha+\tau} \right) \ge A|x|^{-\alpha} $ a.e. in $B_R$. This proves the assertion. \end{proof}

\begin{proof}[Proof of Theorem~\ref{thm:origin}\textup{(b)}] The solutions constructed in Theorems~\ref{thm:main2} and \ref{thm:main3} belong to $H^1_{\rm loc}(\Omega)$. Moreover, since $g\ge0$, their weak formulation implies that $$ \cE(u,\psi) - \mu\int_\Omega\frac{u\psi}{|x|^2}\dd x = \int_\Omega\frac{g\psi}{u^\gamma}\dd x \ge0 $$ for every nonnegative $\psi\in C^\infty_c(\Omega)$. Thus, they satisfy \eqref{eq:supersol}, and part \textup{(a)} applies. Let now $u$ be the solution constructed in Theorem~\ref{thm:main1} as the pointwise limit of the solutions $(u_n)$ of \eqref{problem11}. Set $ \widetilde g:=\min\{g,1\}. $ Then $$ \widetilde g\in L^\infty(\Omega), \qquad 0\le\widetilde g\le g, \qquad \widetilde g\not\equiv0. $$ Let $\widetilde u_n$ be the solution of \eqref{problem11} corresponding to the datum $\widetilde g$. Since $ \widetilde g_n\le g_n, $ the comparison argument used in the proof of Lemma~\ref{lem:monotone}, applied to $(\widetilde u_n-u_n)^+$, gives $ \widetilde u_n\le u_n$ a.e. in $\Omega. $ Passing to the pointwise limit, we obtain $ \widetilde u:=\lim_{n\to\infty}\widetilde u_n \le u $ a.e. in $\Omega. $ Since $\widetilde g\in L^\infty(\Omega)$ and $\mu<\LN$, Theorem~\ref{thm:main2} applies to the datum $\widetilde g$. Hence $\widetilde u\in H^1_0(\Omega)$ is a weak solution of the corresponding problem. By part \textup{(a)}, there exist $R>0$ and $A>0$ such that $ \widetilde u(x)\ge A|x|^{-\alpha_\mu} $ for a.e. $x\in B_R$. Since $u\ge\widetilde u$, the same estimate holds for $u,$ that is, $ u(x)\ge A|x|^{-\alpha_\mu} $ for a.e. $x\in B_R$. Finally, let $r>0$ be such that $B_r\subset\Omega$, and suppose that $p\alpha_\mu\ge N$. Then $$ \begin{aligned} \int_{B_r}u^p\dd x \ge A^p\int_{B_{\min\{r,R\}}}|x|^{-p\alpha_\mu}\dd x = A^p|\mathbb S^{N-1}| \int_0^{\min\{r,R\}} t^{N-1-p\alpha_\mu}\dd t = \infty. \end{aligned} $$ 

Therefore, $ u\notin L^p(B_r)$ whenever $p\ge\frac{N}{\alpha_\mu}. $ In particular, $u$ is unbounded in every neighborhood of the origin. \end{proof}

\begin{proof}[Proof of Theorem~\ref{thm:origin}\textup{(c)}] Under the assumptions of the theorem, the solution $u$ is provided by Theorem~\ref{thm:main2} when $\gamma\le1$, and by Theorem~\ref{thm:main3} when $\gamma>1$. Notice that, when $\gamma\le1$, the condition $\delta_m\ge1$ corresponds to $\gamma=1$ or, if $0\le\gamma<1$, to $m\ge m^\#$. By part \textup{(b)}, there exist $R>0$ and $A>0$ such that $ u(x)\ge A|x|^{-\alpha_\mu} $ for a.e. $x\in B_R$. By \eqref{eq:exponents} and Lemma~\ref{lem:exponents}\textup{(iii)}, $$ \frac{N}{\sigma_m} = \frac{N-2m}{m(1+\gamma)} = \frac{N-2}{1+\delta_m}. $$ Moreover, $$ \begin{aligned} f\left(\frac{N-2}{1+\delta_m}\right) = \frac{N-2}{1+\delta_m} \left( N-2-\frac{N-2}{1+\delta_m} \right) = \frac{(N-2)^2\delta_m}{(1+\delta_m)^2} = \mu_{m,\gamma}. \end{aligned} $$ 

Since $\delta_m\ge1$, we have $ \frac{N-2}{1+\delta_m}\le\frac{N-2}{2}. $ Thus, $\frac{N-2}{1+\delta_m}$ is the smaller root of $ f(\alpha)=\mu_{m,\gamma}, $ and therefore $$ \alpha_{\mu_{m,\gamma}} = \frac{N-2}{1+\delta_m} = \frac{N}{\sigma_m}. $$ Since the map $\mu\mapsto\alpha_\mu$ is increasing and $\mu\ge\mu_{m,\gamma}$, it follows that $ \alpha_\mu \ge \alpha_{\mu_{m,\gamma}} = \frac{N}{\sigma_m}. $ Equivalently, $ \sigma_m\alpha_\mu\ge N. $ Part \textup{(b)}, applied with $p=\sigma_m$, now gives $ u\notin L^{\sigma_m}(B_R). $ Consequently, $ u\notin L^{\sigma_m}(\Omega). $ This proves the optimality assertion. \end{proof}

\subsection{The critical case: proof of Theorem \ref{thm:critical}}
Throughout this subsection $\mu=\LN$. We use the following improved Hardy inequality, due to V\'azquez and Zuazua \cite{vazquez2000hardy} (see also \cite{brezis1997blow} for the corresponding $L^p$ version): for every $q\in[1,2)$ there exists $C_q=C_q(N,\Omega)>0$ such that
\begin{equation}\label{eq:improvedhardy}
\int_\Omega|\nabla u|^2\dd x-\LN\int_\Omega\frac{u^2}{|x|^2}\dd x\ge C_q\,\|\nabla u\|_q^2\qquad\text{for all }u\in H^1_0(\Omega).
\end{equation}

\begin{lemma}\label{lem:critical-coercive} For every $n\in\N$, the constant $$ c_n := \inf_{u\in H^1_0(\Omega)\setminus\{0\}} \frac{ \rho(u)^2-\LN\displaystyle\int_\Omega V_nu^2\dd x }{ \rho(u)^2 } $$ 
is strictly positive, where $V_n(x):=\frac{1}{|x|^2+1/n}.$ Consequently, Lemmas~\ref{lem:linear}, \ref{lem:Pn}, and \ref{lem:monotone}, as well as Remark~\ref{rem:limit} and Proposition~\ref{prop:limit}, remain valid when $\mu=\LN$. In the proof of Lemma~\ref{lem:Pn}, the coercivity constant $1-\mu/\LN$ is replaced by $c_n$, while in the comparison argument of Lemma~\ref{lem:monotone} it is replaced by $c_{n+1}$. \end{lemma}

\begin{proof} Recall that $V_n(x):=\frac{1}{|x|^2+1/n}$ and then by relation \eqref{eq:Vn} and the Hardy inequality \eqref{eq:hardy}, we have $$ \LN\int_\Omega V_nu^2\dd x \le \LN\int_\Omega\frac{u^2}{|x|^2}\dd x \le \|\nabla u\|_2^2 \le \rho(u)^2 $$ for every $u\in H^1_0(\Omega)$. Hence $c_n\ge0$. 

Suppose by contradiction that $c_n=0$. Then there exists a sequence $(u_k)\subset H^1_0(\Omega)$ such that $ \rho(u_k)=1 $ and $$ \rho(u_k)^2 - \LN\int_\Omega V_nu_k^2\dd x \longrightarrow0. $$ Equivalently, $ \LN\int_\Omega V_nu_k^2\dd x \longrightarrow1. $ Since $(u_k)$ is bounded in $H^1_0(\Omega)$, after passing to a subsequence, there exists $u\in H^1_0(\Omega)$ such that $ u_k\rightharpoonup u$ weakly in $H^1_0(\Omega) $ and, by the Rellich--Kondrachov theorem, $ u_k\to u $ strongly in $L^2(\Omega). $ Since $V_n\in L^\infty(\Omega)$, it follows that $$ \LN\int_\Omega V_nu^2\dd x = \lim_{k\to\infty} \LN\int_\Omega V_nu_k^2\dd x = 1. $$ 
In particular, $u\not\equiv0$. By weak lower semicontinuity, $ \rho(u)^2 \le \liminf_{k\to\infty}\rho(u_k)^2 = 1. $ 

On the other hand, $$ \begin{aligned} 1 = \LN\int_\Omega V_nu^2\dd x \le \LN\int_\Omega\frac{u^2}{|x|^2}\dd x \le \|\nabla u\|_2^2 \le \rho(u)^2 \le 1. \end{aligned} $$

Therefore, equality holds throughout. In particular, $$ \|\nabla u\|_2^2 - \LN\int_\Omega\frac{u^2}{|x|^2}\dd x = 0. $$ Applying the improved Hardy inequality \eqref{eq:improvedhardy}, with any fixed $q\in[1,2)$, we obtain $ C_q\|\nabla u\|_q^2\le0. $ Hence $\nabla u=0$ a.e. in $\Omega$. Since $u\in H^1_0(\Omega)$, the Poincar\'e inequality gives $ u=0$ a.e. in $\Omega, $ contradicting $ \LN\int_\Omega V_nu^2\dd x=1. $ Thus, $ c_n>0. $ 

It follows that the bilinear form associated with the $n$-th approximating problem satisfies $$ \cB_n(v,v) = \rho(v)^2 - \LN\int_\Omega V_nv^2\dd x \ge c_n\rho(v)^2 $$ for every $v\in H^1_0(\Omega)$. Therefore, the existence, uniqueness, positivity, and boundedness arguments in Lemmas~\ref{lem:linear} and \ref{lem:Pn} remain valid, with $c_n$ replacing $1-\mu/\LN$ wherever coercivity is used. For the monotonicity argument, \eqref{eq:Vn} gives $$ \begin{aligned} \rho(z^+)^2 - \LN\int_\Omega V_{n+1}(z^+)^2\dd x &\ge c_{n+1}\rho(z^+)^2. \end{aligned} $$ 

Hence the proof of Lemma~\ref{lem:monotone} remains valid with $c_{n+1}$ in place of $1-\mu/\LN$. Remark~\ref{rem:limit} then follows without change. Finally, Proposition~\ref{prop:limit} only uses the monotonicity and positivity properties of $(u_n)$, the assumed compactness estimates, and Lemma~\ref{lem:welldef}; thus its conclusion also remains valid in the critical case. \end{proof}

\begin{lemma}\label{lem:critical-estimate} Under the assumptions of Theorem~\ref{thm:critical}, let $(u_n)$ be the sequence of solutions of \eqref{problem11} corresponding to $\mu=\LN$. Then, for every $q\in[1,2)$, $(u_n)$ is bounded in $W^{1,q}_0(\Omega)$. Moreover, the sequence $([u_n]_s)$ is bounded. \end{lemma} 

\begin{proof} Testing \eqref{weak11} with $u_n$ and using $ V_n\le\frac{1}{|x|^2}, $ we obtain $$ \begin{aligned} \rho(u_n)^2 -\LN\int_\Omega V_nu_n^2\dd x = \int_\Omega \frac{g_nu_n}{(u_n+1/n)^\gamma}\dd x. \end{aligned} $$ 

Moreover, for every $q\in[1,2)$, the improved Hardy inequality \eqref{eq:improvedhardy} gives $$ \begin{aligned} \rho(u_n)^2 -\LN\int_\Omega V_nu_n^2\dd x &\ge \|\nabla u_n\|_2^2 -\LN\int_\Omega\frac{u_n^2}{|x|^2}\dd x +\frac{\Cs}{2}[u_n]_s^2 \\ &\ge C_q\|\nabla u_n\|_q^2 +\frac{\Cs}{2}[u_n]_s^2. \end{aligned} $$ Consequently, \begin{equation}\label{eq:crit-energy} C_q\|\nabla u_n\|_q^2 +\frac{\Cs}{2}[u_n]_s^2 \le \int_\Omega \frac{g_nu_n}{(u_n+1/n)^\gamma}\dd x. \end{equation} Suppose first that $\gamma=1$. Since $$ 0\le\frac{u_n}{u_n+1/n}\le1, $$ we have $$ \int_\Omega \frac{g_nu_n}{u_n+1/n}\dd x \le \int_\Omega g_n\dd x \le \|g\|_1. $$ 

It follows from \eqref{eq:crit-energy} that $ \|\nabla u_n\|_q\le C_q $ for every $q\in[1,2)$, and $ [u_n]_s\le C. $ The Poincar\'e inequality then yields the boundedness of $(u_n)$ in $W^{1,q}_0(\Omega)$. Assume now that $0\le\gamma<1$ and $m>m^\#$. By the definition of $m^\#$, $ (1-\gamma)m'<2^*. $ Since $$ q^*:=\frac{Nq}{N-q}\longrightarrow2^* \qquad\text{as }q\to2^-, $$ we may choose $q_0\in(1,2)$ such that $ q_0^*\ge(1-\gamma)m'. $ By H\"older's inequality, $$ \begin{aligned} \int_\Omega \frac{g_nu_n}{(u_n+1/n)^\gamma}\dd x \le \int_\Omega gu_n^{1-\gamma}\dd x \le \|g\|_m \|u_n\|_{(1-\gamma)m'}^{1-\gamma}. \end{aligned} $$ 

Since $\Omega$ is bounded and $q_0^*\ge(1-\gamma)m'$, the Sobolev and Poincar\'e inequalities give $$ \|u_n\|_{(1-\gamma)m'} \le C\|u_n\|_{q_0^*} \le C\|\nabla u_n\|_{q_0}. $$ Therefore, $$ \int_\Omega \frac{g_nu_n}{(u_n+1/n)^\gamma}\dd x \le C\|g\|_m \|\nabla u_n\|_{q_0}^{1-\gamma}. $$ Taking $q=q_0$ in \eqref{eq:crit-energy} and discarding the nonnegative fractional term, we obtain $$ C_{q_0}\|\nabla u_n\|_{q_0}^2 \le C\|g\|_m \|\nabla u_n\|_{q_0}^{1-\gamma}. $$ Hence, either $\|\nabla u_n\|_{q_0}=0$, or $ \|\nabla u_n\|_{q_0}^{1+\gamma} \le C\|g\|_m. $ Thus, $ \|\nabla u_n\|_{q_0}\le C. $ It follows that the right-hand side of \eqref{eq:crit-energy} is bounded independently of $n$. Returning to \eqref{eq:crit-energy}, we obtain, for every fixed $q\in[1,2)$, $$ \|\nabla u_n\|_q\le C_q \qquad\text{and}\qquad [u_n]_s\le C. $$ Finally, the Poincar\'e inequality gives $ \|u_n\|_{W^{1,q}_0(\Omega)}\le C_q $ for every $q\in[1,2)$. This completes the proof. \end{proof}

\begin{lemma}\label{lem:monotone-mu} Let $(m,\gamma)$ and $g$ be as in Theorem~\ref{thm:critical}, let $0<\mu<\LN$, and let $u_{n,\mu}$ and $u_{n,\LN}$ denote the solutions of \eqref{problem11} corresponding to the parameters $\mu$ and $\LN$, respectively. Then, for every $n\in\N$, $$ u_{n,\mu}\le u_{n,\LN} \quad\text{a.e. in }\Omega. $$ \end{lemma}

\begin{proof} Fix $n\in\N$ and set $ z:=u_{n,\mu}-u_{n,\LN}. $ Subtracting the weak formulations satisfied by $u_{n,\mu}$ and $u_{n,\LN}$ and testing the resulting identity with $z^+\in H^1_0(\Omega)$, we obtain $$ \begin{aligned} \cE(z,z^+) -\LN\int_\Omega V_nzz^+\dd x &= (\mu-\LN) \int_\Omega V_nu_{n,\mu}z^+\dd x \\ &\qquad\quad+ \int_\Omega g_n \left( \frac{1}{(u_{n,\mu}+1/n)^\gamma} - \frac{1}{(u_{n,\LN}+1/n)^\gamma} \right)z^+\dd x. \end{aligned} $$ Since $\mu-\LN<0$ and $ V_nu_{n,\mu}z^+\ge0, $ the first term on the right-hand side is nonpositive. Moreover, on the set $ \{z>0\} = \{u_{n,\mu}>u_{n,\LN}\}, $ the monotonicity of the function $ t\longmapsto\frac{1}{(t+1/n)^\gamma} $ gives $$ \frac{1}{(u_{n,\mu}+1/n)^\gamma} - \frac{1}{(u_{n,\LN}+1/n)^\gamma} \le0. $$ Thus, the second term on the right-hand side is also nonpositive, and hence $$ \cE(z,z^+) -\LN\int_\Omega V_nzz^+\dd x \le0. $$ On the other hand, Lemma~\ref{lem:alg}\textup{(iii)} gives $ \cE(z,z^+)\ge\rho(z^+)^2. $ Since $ zz^+=(z^+)^2$ a.e. in $\Omega, $ Lemma~\ref{lem:critical-coercive} yields $$ \begin{aligned} \cE(z,z^+) -\LN\int_\Omega V_nzz^+\dd x &\ge \rho(z^+)^2 -\LN\int_\Omega V_n(z^+)^2\dd x \ge c_n\rho(z^+)^2. \end{aligned} $$

Consequently, $ c_n\rho(z^+)^2\le0. $ Since $c_n>0$, it follows that $ z^+=0$ a.e. in $\Omega. $ Therefore, $ u_{n,\mu}\le u_{n,\LN} $ a.e. in $\Omega. $ \end{proof}

\begin{proof}[Proof of Theorem~\ref{thm:critical}] By Lemma~\ref{lem:critical-estimate}, $(u_n)$ is bounded in $W^{1,q}_0(\Omega)$ for every $q\in[1,2)$. Fix $ q\in\left(\frac{N}{N-1},2\right). $ Since $(u_n)$ is also bounded in $L^1(\Omega)$, Proposition~\ref{prop:limit}, which remains valid for $\mu=\LN$ by Lemma~\ref{lem:critical-coercive}, shows that $ u:=\lim_{n\to\infty}u_n $ is a weak solution of \eqref{problem1}. Let now $q\in(1,2)$ be arbitrary. Since $(u_n)$ is bounded in the reflexive space $W^{1,q}_0(\Omega)$, every subsequence admits a further subsequence converging weakly in $W^{1,q}_0(\Omega)$.

On the other hand, $ u_n\to u$ in $L^1(\Omega). $ Therefore, every weak limit coincides with $u$, and hence $ u\in W^{1,q}_0(\Omega) $ and $ u_n\rightharpoonup u$ weakly in $W^{1,q}_0(\Omega). $ Since this holds for every $q\in(1,2)$ and $\Omega$ is bounded, we also have $ u\in W^{1,1}_0(\Omega). $ Thus, $ u\in W^{1,q}_0(\Omega)$ for every $q\in[1,2). $ For every $p<2^*$, we may choose $q\in(1,2)$ sufficiently close to $2$ so that $ p\le q^*:=\frac{Nq}{N-q}. $ The Sobolev embedding then gives $ u\in L^p(\Omega). $ Moreover, Lemma~\ref{lem:critical-estimate} gives $ [u_n]_s\le C. $ Since $u_n\to u$ a.e. in $\Omega$, we also have $ u_n(x)-u_n(y)\to u(x)-u(y) $ for a.e. $(x,y)\in\R^{2N}$. Fatou's lemma therefore yields $$ \begin{aligned} [u]_s^2 = \iint_{\R^{2N}} \frac{|u(x)-u(y)|^2}{|x-y|^{N+2s}}\dd x\dd y \le \liminf_{n\to\infty} \iint_{\R^{2N}} \frac{|u_n(x)-u_n(y)|^2}{|x-y|^{N+2s}}\dd x\dd y = \liminf_{n\to\infty}[u_n]_s^2 <\infty. \end{aligned} $$
This proves \textup{(i)}.

We now prove \textup{(ii)}. Let $ \beta<\frac{N-2}{2}. $ Choose $\lambda\in(0,\LN)$ sufficiently close to $\LN$ so that $ \alpha_\lambda\ge\beta. $ This is possible because $ \alpha_\lambda\longrightarrow\frac{N-2}{2}$ as $\lambda\to\LN^-. $ Let $u_\lambda$ be the solution constructed in Theorem~\ref{thm:main2} for the parameter $\lambda$ and the same datum $g$. Such a solution exists because the assumptions on $(m,\gamma)$ in Theorem~\ref{thm:critical} imply the assumptions of Theorem~\ref{thm:main2}. By Lemma~\ref{lem:monotone-mu}, $ u_{n,\lambda}\le u_{n,\LN}$ a.e. in $\Omega $ for every $n\in\N$. Passing to the pointwise limit gives $ u_\lambda\le u$ a.e. in $\Omega. $ By Theorem~\ref{thm:origin}\textup{(b)}, there exist $R>0$ and $A>0$ such that $ u_\lambda(x)\ge A|x|^{-\alpha_\lambda} $ for a.e. $x\in B_R$. Decreasing $R$ if necessary, we may assume that $R<1$. Since $\alpha_\lambda\ge\beta$, we have $ |x|^{-\alpha_\lambda}\ge|x|^{-\beta} $ for $0<|x|<R$. Consequently, $$ u(x)\ge u_\lambda(x)\ge A|x|^{-\beta} $$ for a.e. $x\in B_R$. Finally, let $p>2^*$. Then $ \frac{N}{p}<\frac{N-2}{2}. $ Applying the preceding estimate with $ \beta:=\frac{N}{p}, $ we obtain constants $R>0$ and $A>0$ such that $ u(x)\ge A|x|^{-N/p} $ for a.e. $x\in B_R$. Hence, for every $r>0$ such that $B_r\subset\Omega$, $$ \begin{aligned} \int_{B_r}u^p\dd x \ge A^p \int_{B_{\min\{r,R\}}}|x|^{-N}\dd x = A^p|\mathbb S^{N-1}| \int_0^{\min\{r,R\}}\frac{\dd t}{t} = \infty. \end{aligned} $$

Therefore, $ u\notin L^p(B_r) $ for every $p>2^*$ and every $r>0$ such that $B_r\subset\Omega$. This proves \textup{(ii)} and hence the theorem. \end{proof}

\begin{remark}\label{rem:critical} In the purely local case, solutions of $$ -\Delta u-\LN\frac{u}{|x|^2}=f, \qquad f\ge0,\quad f\not\equiv0, $$ need not belong to $H^1_0(\Omega)$; see, for example, \cite{peral2021elliptic,vazquez2000hardy}. This explains why, when $\mu=\LN$, the natural energy space is the completion of $C^\infty_c(\Omega)$ with respect to the quadratic-form norm $$ \left( \rho(u)^2 - \LN\int_\Omega\frac{u^2}{|x|^2}\dd x \right)^{1/2}. $$ The approximate solutions $(u_n)$ are bounded in this space by \eqref{eq:crit-energy}. A more precise description of the behavior of the solution near the origin in the critical case for the mixed operator remains an interesting open problem. \end{remark}

\section{Appendix}
In this section, we give the proofs of Lemmas \ref{lem:alg}, \ref{lem:stampacchia} and \ref{lem:strongmin} to make the paper self-contained.

\begin{proof}[Proof of Lemma \ref{lem:alg}] \textup{(i)} The inequality is well known; see, for example, \cite[Lemma A.4]{arora2023combined}. We include its short proof for completeness. By symmetry, it suffices to consider $a>b\ge 0$. The Cauchy--Schwarz inequality gives

\begin{align*} \left( a^{\frac{\theta+1}{2}}-b^{\frac{\theta+1}{2}} \right)^2 &= \left( \frac{\theta+1}{2} \int_b^a t^{\frac{\theta-1}{2}}\dd t \right)^2 \le \frac{(\theta+1)^2}{4} (a-b)\int_b^a t^{\theta-1}\dd t \\ &= \frac{(\theta+1)^2}{4\theta} (a-b)(a^\theta-b^\theta), \end{align*} which proves \eqref{eq:alg}. 

\textup{(ii)} Since $\Phi$ is nondecreasing, $a-b$ and $\Phi(a)-\Phi(b)$ have the same sign, and therefore \[ (a-b)(\Phi(a)-\Phi(b))\ge 0. \] 

Suppose now that $\Phi$ is $1$-Lipschitz. If $a=b$, the second inequality is immediate. If $a\ne b$, then \[ 0\le \frac{\Phi(a)-\Phi(b)}{a-b} \le 1. \] 

Using $t^2\le t$ for $t\in[0,1]$ and multiplying by $(a-b)^2$, we obtain 
\[ \left(\Phi(a)-\Phi(b)\right)^2 \le (a-b)(\Phi(a)-\Phi(b)). \] 

\textup{(iii)} The standard truncation properties of Sobolev functions yield $z^\pm\in H^1_0(\Omega)$ and \[ \nabla z^+ = \chi_{\{z>0\}}\nabla z, \qquad \nabla z^- = -\chi_{\{z<0\}}\nabla z \quad\text{a.e. in }\Omega. \] Consequently, $\nabla z\cdot\nabla z^+ = |\nabla z^+|^2$ and $\nabla z\cdot\nabla(-z^-) = |\nabla z^-|^2$ a.e. in $\Omega.$ For the nonlocal part, using $z=z^+-z^-$ and $z^+z^-=0$, we obtain, for every $x,y\in\RN$, 
\begin{align*} (z(x)-z(y))(z^+(x)-z^+(y))  &= (z^+(x)-z^+(y))^2 +z^-(x)z^+(y)+z^-(y)z^+(x) \\ &\quad\ge (z^+(x)-z^+(y))^2, \end{align*}
because $z^\pm\ge 0$. Similarly, 
\begin{align*} -(z(x)-z(y))(z^-(x)-z^-(y)) & = (z^-(x)-z^-(y))^2 +z^+(x)z^-(y)+z^+(y)z^-(x) \\ &\quad\ge (z^-(x)-z^-(y))^2. \end{align*} Combining these pointwise inequalities with the corresponding local ones gives 
$$\cE(z,z^+)\ge \rho(z^+)^2, \qquad \cE(z,-z^-)\ge \rho(z^-)^2. $$ \end{proof}

\begin{proof}[Proof of Lemma \ref{lem:stampacchia}]
\textup{(a)}  For $k\ge 0$, set $ A_k:=\{x\in\Omega:w(x)>k\} $ and let $G_k(t):=(t-k)^+$. Then $G_k(w)\in H^1_0(\Omega)$ and $$ \nabla w\cdot\nabla G_k(w) = |\nabla G_k(w)|^2 \quad\text{a.e. in }\Omega. $$ Moreover, since $G_k$ is nondecreasing, Lemma~\ref{lem:alg}\textup{(ii)} shows that the nonlocal part of $\cE(w,G_k(w))$ is nonnegative. Testing the equation with $G_k(w)$ and using H\"older's inequality, we obtain 
$$ \begin{aligned} \|\nabla G_k(w)\|_2^2 \le \int_{A_k}fG_k(w)\dd x \le \|f\|_p\|G_k(w)\|_{p'} \le \|f\|_p\|G_k(w)\|_{2^*} |A_k|^{1-\frac1p-\frac1{2^*}}. \end{aligned} $$

By the Sobolev inequality \eqref{eq:sobolev}, $$ \|G_k(w)\|_{2^*}^2 \le S_N\|\nabla G_k(w)\|_2^2. $$ Therefore, $$ \|G_k(w)\|_{2^*} \le S_N\|f\|_p |A_k|^{1-\frac1p-\frac1{2^*}}. $$ For $h>k$, we have $G_k(w)\ge h-k$ on $A_h$, and hence $ (h-k)|A_h|^{1/2^*} \le \|G_k(w)\|_{2^*}. $ Consequently, $$ |A_h| \le \frac{(S_N\|f\|_p)^{2^*}}{(h-k)^{2^*}} |A_k|^\beta, \qquad \beta:= 2^*\left(1-\frac1p-\frac1{2^*}\right). $$ 

Since $ \beta>1 \quad\Longleftrightarrow\quad p>\frac N2, $ the Stampacchia iteration lemma \cite[Lemma B.1, Chapter II]{kinderlehrer1980introduction} yields $|A_d|=0$, where $ d:= 2^{\frac{\beta}{\beta-1}} S_N\|f\|_p |\Omega|^{\frac{\beta-1}{2^*}}. $ Thus $0\le w\le d$ a.e. in $\Omega$, and therefore $$ \|w\|_\infty \le 2^{\frac{\beta}{\beta-1}} S_N|\Omega|^{\frac{\beta-1}{2^*}}\|f\|_p. $$ 

\textup{(b)} Set $ \theta:=\frac{p'}{2p'-2^*}. $ Since $(2^*)'\le p<N/2$, we have $ \frac{N}{N-2}<p'\le 2^*. $ In particular, $2p'>2^*$ and $\theta\ge 1$. Moreover, $ (2\theta-1)p'=\theta2^*. $

For $k>0$, let $ T_k(t):=\min\{t,k\},$ for $t\ge 0$ and define $ \Phi_k(t):=T_k(t)^{2\theta-1}. $ The function $\Phi_k$ is nondecreasing and Lipschitz on $[0,\infty)$, with $\Phi_k(0)=0$. Consequently, $\Phi_k(w)\in H^1_0(\Omega)$ and $$ \begin{aligned} \nabla w\cdot\nabla\Phi_k(w) = (2\theta-1)T_k(w)^{2\theta-2} |\nabla T_k(w)|^2 = \frac{2\theta-1}{\theta^2} \left|\nabla\bigl(T_k(w)^\theta\bigr)\right|^2 \end{aligned} $$ 
a.e. in $\Omega$. Testing the equation with $\Phi_k(w)$ and using Lemma~\ref{lem:alg}\textup{(ii)}, we may discard the nonnegative nonlocal contribution. Thus, $$ \frac{2\theta-1}{\theta^2} \int_\Omega \left|\nabla\bigl(T_k(w)^\theta\bigr)\right|^2\dd x \le \int_\Omega fT_k(w)^{2\theta-1}\dd x. $$ 

By H\"older's inequality and the identity $(2\theta-1)p'=\theta2^*$, we have $$ \begin{aligned} \int_\Omega fT_k(w)^{2\theta-1}\dd x \le \|f\|_p \left( \int_\Omega T_k(w)^{(2\theta-1)p'}\dd x \right)^{1/p'} = \|f\|_p \left( \int_\Omega T_k(w)^{\theta2^*}\dd x \right)^{1/p'}. \end{aligned} $$

Applying the Sobolev inequality to $T_k(w)^\theta$, we obtain $$ \frac{2\theta-1}{\theta^2S_N} \left( \int_\Omega T_k(w)^{\theta2^*}\dd x \right)^{2/2^*} \le \|f\|_p \left( \int_\Omega T_k(w)^{\theta2^*}\dd x \right)^{1/p'}. $$ 

Since $p<N/2$, $ \frac{2}{2^*}-\frac1{p'} = \frac1p-\frac2N >0. $ It follows that $$ \left( \int_\Omega T_k(w)^{\theta2^*}\dd x \right)^{ \frac{2}{2^*}-\frac1{p'} } \le \frac{\theta^2S_N}{2\theta-1}\|f\|_p. $$

A direct computation gives $ \theta2^* = \left(\frac1p-\frac2N\right)^{-1} = \frac{Np}{N-2p} = p^{**}. $ Hence $ \frac{2}{2^*}-\frac1{p'} = \frac1{p^{**}}, $ and therefore $$ \|T_k(w)\|_{p^{**}} \le \frac{\theta^2S_N}{2\theta-1}\|f\|_p. $$ 

The right-hand side is independent of $k$. Since $T_k(w)\uparrow w$ pointwise as $k\to\infty$, the monotone convergence theorem yields $ w\in L^{p^{**}}(\Omega) $ and $$ \|w\|_{p^{**}} \le \frac{\theta^2S_N}{2\theta-1}\|f\|_p. $$ This completes the proof. \end{proof}
 The positivity of the approximate solutions relies on the weak Harnack inequality of Garain and Kinnunen \cite{garain2022regularity}.

\begin{proof}[Proof of Lemma \ref{lem:strongmin}] We use the weak Harnack inequality for nonnegative weak supersolutions of mixed local--nonlocal equations established in \cite{garain2022regularity}. Since $w\ge 0$ in $\RN$, the nonlocal tail of the negative part of $w$ vanishes. Consequently, there exist $t>0$ and $C\ge 1$, possibly depending on $N$, $s$, the coefficient of the nonlocal term, and the radius $r$, such that \begin{equation}\label{eq:weakharnack} \left( \frac{1}{|B_r(x_0)|} \int_{B_r(x_0)}w^t\dd x \right)^{1/t} \le C\essinf_{B_r(x_0)}w \end{equation} whenever $B_{4r}(x_0)\Subset\Omega$. The precise ratio between the radii and the dependence of $C$ on $r$ are irrelevant for the argument. We only use the fact that the left-hand side of \eqref{eq:weakharnack} is positive whenever $w$ does not vanish a.e. in $B_r(x_0)$. Define $$ A:= \left\{ x\in\Omega: \begin{array}{l} \text{there exists }r>0\text{ such that } B_{4r}(x)\Subset\Omega\\ \text{and }\essinf_{B_r(x)}w>0 \end{array} \right\}. $$ 

By definition, $A$ is open. We claim that $A$ is nonempty. Indeed, since $w\not\equiv 0$ in $\Omega$, there exists a ball $B_r(x_0)$ such that $ B_{4r}(x_0)\Subset\Omega $ and $w$ does not vanish a.e. in $B_r(x_0)$. Therefore, $$ \left( \frac{1}{|B_r(x_0)|} \int_{B_r(x_0)}w^t\dd x \right)^{1/t}>0. $$ 

It follows from \eqref{eq:weakharnack} that $ \essinf_{B_r(x_0)}w>0, $ and hence $x_0\in A$. 

We next show that $A$ is relatively closed in $\Omega$. Let $(x_j)\subset A$ be such that $x_j\to x\in\Omega$. Choose $r>0$ such that $ B_{4r}(x)\Subset\Omega. $ For all sufficiently large $j$, we have $x_j\in B_r(x)$. Since $x_j\in A$, there exist $r_j>0$ and $c_j>0$ such that $ w\ge c_j$ a.e. in  $B_{r_j}(x_j). $ The set $ B_{r_j}(x_j)\cap B_r(x) $ has positive measure. Therefore, 
$$ \int_{B_r(x)}w^t\dd x>0. $$ 

Applying \eqref{eq:weakharnack} to the ball $B_r(x)$, we obtain $ \essinf_{B_r(x)}w>0. $ Thus $x\in A$, and $A$ is relatively closed in $\Omega$. Since $\Omega$ is connected and $A$ is nonempty, open, and relatively closed in $\Omega$, it follows that $ A=\Omega. $ Now let $\omega\Subset\Omega$. For every $x\in\overline{\omega}$, there exist $r_x>0$ and $c_x>0$ such that $ w\ge c_x$ a.e. in  $B_{r_x}(x). $ By compactness, there exist $x_1,\ldots,x_m\in\overline{\omega}$ such that $$ \overline{\omega} \subset \bigcup_{i=1}^m B_{r_{x_i}}(x_i). $$ Setting $ c_\omega:=\min_{1\le i\le m}c_{x_i}>0, $ we conclude that $ w\ge c_\omega$ a.e. in $\omega. $ \end{proof}

\section*{Acknowledgements}
SG gratefully acknowledges the financial support for this research work under ARG-MATRICS, grant No: ANRF/ARGM/2025/001570/MTR, Anusandhan National Research Foundation (ANRF), Government of India. After the original submission of the present manuscript to a journal, we became aware of the interesting preprint by Abdellaoui, Biroud, Laamri and Montoro \cite{ABLM26}, which independently addresses, among other problems, the singular problem considered here. The present work develops a different analytical framework, with particular emphasis on quantitative regularity, explicit summability exponents and the optimality of the associated Hardy thresholds.

\section*{Declarations}
\noindent\textbf{Data availability.} Data sharing is not applicable to this article, as no datasets were generated or analyzed.

\noindent\textbf{Conflict of interest.} The authors declare that they have no conflict of interest.

\noindent\textbf{AI Declaration}
The authors used ChatGPT (OpenAI, GPT-6) partially for minor language editing in a few places and take full responsibility for the reviewed manuscript.

\end{document}